\documentclass{article}
\usepackage[utf8]{inputenc}
\usepackage[margin=1in]{geometry}
\usepackage{todonotes}
\usepackage{amsmath}
\usepackage{amsthm}
\usepackage{amssymb}
\usepackage{hyperref}
\usepackage{caption}
\usepackage{subcaption}
\usepackage[title]{appendix}
\usepackage{multirow}
\usepackage{commath}
\usepackage{mathtools}
\usepackage{comment}
\usepackage{pdflscape}
\usepackage{soul}
\usepackage{adjustbox}
\usepackage{enumitem}
\usepackage{xcolor}
\usepackage{pifont}
\usepackage[linesnumbered, ruled]{algorithm2e}
\usepackage{todonotes}
\usepackage{tikz}
\usepackage{tikz-3dplot}
\usepackage{pgfplots}

\theoremstyle{definition}

\numberwithin{equation}{section}

\title{Generation of Paths for Motion Planning for a Dubins Vehicle on Sphere\thanks{DISTRIBUTION STATEMENT A. Approved for public release. Distribution is unlimited. AFRL-2025-0643; Cleared 05 Feb 2025.}}

\author{Deepak Prakash Kumar\thanks{Deepak Prakash Kumar and Swaroop Darbha are with the Department of Mechanical Engineering, Texas A\&M University, College Station, TX 77843, USA (e-mail: {\tt\footnotesize deepakprakash1997@gmail.com, dswaroop@tamu.edu}).}, Swaroop Darbh$\text{a}^\dagger$, Satyanarayana Gupta Manyam\thanks{Satyanarayana Gupta Manyam is with the Infoscitex Corp., 4027 Col Glenn Hwy, Dayton, OH 45431, USA (e-mail: {\tt\footnotesize msngupta@gmail.com})}, David Casbeer\thanks{David Casbeer is with the Autonomous Control Branch, Air Force Research Laboratory, Wright-Patterson Air Force Base, OH 45433 USA (e-mail:
{\tt\footnotesize david.casbeer@afresearchlab.com})}
}

\date{\today}

\begin{document}

\maketitle

\begin{abstract}
    In this article, the candidate optimal paths for a Dubins vehicle on a sphere are analytically derived. In particular, the arc angles for segments in $CGC$, $CCC$, $CCCC$, and $CCCCC$ paths, which have previously been shown to be optimal depending on the turning radius $r$ of the vehicle in \cite{kumar2025newapproachmotionplanning}, are analytically derived. The derived expressions are used for the implementation provided in https://github.com/DeepakPrakashKumar/Motion-planning-on-sphere.
\end{abstract}

\section{Derivation of Closed-Form Expressions for Paths on a Sphere}

Consider the differential equations corresponding to the Sabban frame. As $u_g$ is piecewise constant, the differential equations for each segment can be represented as
\begin{align}
    \begin{pmatrix}
        \mathbf{X}' (s) & \mathbf{T}' (s) & \mathbf{N}' (s)
    \end{pmatrix} &= \begin{pmatrix}
        \mathbf{X} (s) & \mathbf{T} (s) & \mathbf{N} (s)
    \end{pmatrix} \underbrace{\begin{pmatrix}
        0 & -1 & 0 \\
        1 & 0 & -u_g \\
        0 & u_g & 0
    \end{pmatrix}}_{\Omega}.
\end{align}
The solution to the above differential equations can be obtained as
\begin{align}
    \begin{pmatrix}
        \mathbf{X} (s) & \mathbf{T} (s) & \mathbf{N} (s)
    \end{pmatrix} = \begin{pmatrix}
        \mathbf{X} (s_i) & \mathbf{T} (s_i) & \mathbf{N} (s_i)
    \end{pmatrix} \left(e^{\Omega^T \Delta s} \right)^T,
\end{align}
where $\Delta s = s - s_i.$ The expression for $e^{\Omega \Delta s}$ can be obtained using the Euler-Rodriguez formula for the exponential of a skew-symmetric matrix. Moreover, since $s = \phi_G$ for a great circle turn, and $s = r \phi_L$ and $s = r \phi_R$ for the left and right tight turns, respectively, the solution to the differential equations can be written as
\begin{align}
    R_{after} = R_{before} R_{seg},
\end{align}
where $R_{after}$ and $R_{before}$ denote the configurations after and before the segment, respectively. Moreover, $R_{seg}$ represents the rotation matrix corresponding to a chosen segment, and is given by
\begin{align} \label{eq: R_segments}
    R_{seg} = \begin{cases}
        R_G (\phi) = \begin{pmatrix}
        c \phi & - s \phi & 0 \\
        s \phi & c \phi & 0 \\
        0 & 0 & 1
    \end{pmatrix}, & u_g = 0 \\
        R_L (r, \phi) = \begin{pmatrix}
            1 - (1 - c \phi) r^2 & - r s \phi & (1 - c \phi) r \sqrt{1 - r^2} \\
            r s \phi & c \phi & - s \phi \sqrt{1 - r^2} \\
            (1 - c \phi) r \sqrt{1 - r^2} & s \phi \sqrt{1 - r^2} & c \phi + (1 - c \phi) r^2
        \end{pmatrix}, & u_g = U_{max} \\
        R_R (r, \phi) = \begin{pmatrix}
            1 - (1 - c \phi) r^2 & - r s \phi & -(1 - c \phi) r \sqrt{1 - r^2} \\
            r s \phi & c \phi & s \phi \sqrt{1 - r^2} \\
            -(1 - c \phi) r \sqrt{1 - r^2} & -s \phi \sqrt{1 - r^2} & c \phi + (1 - c \phi) r^2
        \end{pmatrix}, & u_g = -U_{max}
    \end{cases}.
\end{align}

For deriving the closed-form expressions for the paths, the following two vectors will be utilized:
\begin{align*}
    \mathbf{u}_{L} := \begin{pmatrix}
        \sqrt{1 - r^2} \\
        0 \\
        r
    \end{pmatrix}, \quad \mathbf{u}_{R} := \begin{pmatrix}
        -\sqrt{1 - r^2} \\
        0 \\
        r
    \end{pmatrix}.
\end{align*}
It should be noted that the vector $\mathbf{u}_{L}$ corresponds to the axial vector of the rotation matrix $R_L (r, \phi).$ Hence, $R_L (r, \phi) \mathbf{u}_{L} = \mathbf{u}_{L}$. Similarly, the vector $\mathbf{u}_{R}$ corresponds to the axial vector of the rotation matrix $R_R (r, \phi).$ Hence, $R_R (r, \phi) \mathbf{u}_{R} = \mathbf{u}_{R}$.

\subsection{LGL path}

The equation to be solved is given by
\begin{align} \label{eq: LGL_general_path_solve}
    R_L (r, \phi_1) R_G (\phi_2) R_L (r, \phi_3) = \begin{pmatrix}
        \alpha_{11} & \alpha_{12} & \alpha_{13} \\
        \alpha_{21} & \alpha_{22} & \alpha_{23} \\
        \alpha_{31} & \alpha_{32} & \alpha_{33}
    \end{pmatrix}.
\end{align}
Pre-multiplying Eq.~\eqref{eq: LGL_general_path_solve} by $\mathbf{u}^T_{L}$ and post-multiplying by $\mathbf{u}_{L},$
\begin{align*}
    \mathbf{u}^T_{L} R_L (r, \phi_1) R_G (\phi_2) R_L (r, \phi_3) \mathbf{u}_{L} = \mathbf{u}^T_{L} R_G (\phi_2) \mathbf{u}_{L} &= \mathbf{u}^T_{L} \begin{pmatrix}
        \alpha_{11} & \alpha_{12} & \alpha_{13} \\
        \alpha_{21} & \alpha_{22} & \alpha_{23} \\
        \alpha_{31} & \alpha_{32} & \alpha_{33}
    \end{pmatrix} \mathbf{u}_{L}. \\
    \therefore \begin{pmatrix}
        \sqrt{1 - r^2} & 0 & r
    \end{pmatrix} \begin{pmatrix}
        c \phi_2 & - s \phi_2 & 0 \\
        s \phi_2 & c \phi_2 & 0 \\
        0 & 0 & 1
    \end{pmatrix} \begin{pmatrix}
        \sqrt{1 - r^2} \\
        0 \\
        r
    \end{pmatrix} &= \begin{pmatrix}
        \sqrt{1 - r^2} & 0 & r
    \end{pmatrix} \begin{pmatrix}
        \alpha_{11} & \alpha_{12} & \alpha_{13} \\
        \alpha_{21} & \alpha_{22} & \alpha_{23} \\
        \alpha_{31} & \alpha_{32} & \alpha_{33}
    \end{pmatrix} \begin{pmatrix}
        \sqrt{1 - r^2} \\
        0 \\
        r
    \end{pmatrix}.
\end{align*}
Simplifying the above equation,
\begin{align} \label{eq: cos_phi_2_LGL_path_expression}
\begin{split}
    (1 - r^2) c \phi_2 + r^2 &= (1 - r^2) \alpha_{11} + r \sqrt{1 - r^2} (\alpha_{13} + \alpha_{31}) + r^2 \alpha_{33}. \\
    \implies c \phi_2 &= \frac{\alpha_{11} + r \sqrt{1 - r^2} (\alpha_{13} + \alpha_{31}) + r^2 (\alpha_{33} - \alpha_{11} - 1)}{1 - r^2},
\end{split}
\end{align}
which yields at most two solutions for $\phi_2 \in [0, 2 \pi)$ if the absolute value of the RHS is less than or equal to 1.

Consider pre-multiplying Eq.~\eqref{eq: LGL_general_path_solve} by $\mathbf{u}_{R}^T$ and post-multiplying by $\mathbf{u}_{L}$, which yields
\begin{align} \label{eq: uRTLGLuL}
    \mathbf{u}^T_{R} R_L (r, \phi_1) R_G (\phi_2) R_L (r, \phi_3) \mathbf{u}_{L} = \mathbf{u}^T_{R} R_L (r, \phi_1) R_G (\phi_2) \mathbf{u}_{L} = \mathbf{u}^T_{R} \begin{pmatrix}
        \alpha_{11} & \alpha_{12} & \alpha_{13} \\
        \alpha_{21} & \alpha_{22} & \alpha_{23} \\
        \alpha_{31} & \alpha_{32} & \alpha_{33}
    \end{pmatrix} \mathbf{u}_{L}.
\end{align}
The LHS of the above equation can be expanded as
\begin{align*}
    &\mathbf{u}^T_{R} R_L (r, \phi_1) R_G (\phi_2) \mathbf{u}_{L} \\
    &= \mathbf{u}^T_{R} R_L (r, \phi_1) \begin{pmatrix}
        c \phi_2 & - s \phi_2 & 0 \\
        s \phi_2 & c \phi_2 & 0 \\
        0 & 0 & 1
    \end{pmatrix} \begin{pmatrix}
        \sqrt{1 - r^2} \\
        0 \\
        r
    \end{pmatrix} \\
    &= \begin{pmatrix}
        -\sqrt{1 - r^2} & 0 & r
    \end{pmatrix} \begin{pmatrix}
        1 - (1 - c \phi_1) r^2 & - r s \phi_1 & (1 - c \phi_1) r \sqrt{1 - r^2} \\
        r s \phi_1 & c \phi_1 & - s \phi_1 \sqrt{1 - r^2} \\
        (1 - c \phi_1) r \sqrt{1 - r^2} & s \phi_1 \sqrt{1 - r^2} & c \phi_1 + (1 - c \phi_1) r^2
    \end{pmatrix} \begin{pmatrix}
        \sqrt{1 - r^2} c \phi_2 \\
        \sqrt{1 - r^2} s \phi_2 \\
        r
    \end{pmatrix} \\
    &= \begin{pmatrix}
        -\sqrt{1 - r^2} + 2 (1 - c \phi_1) r^2 \sqrt{1 - r^2} & 2 r \sqrt{1 - r^2} s \phi_1 & 2 (1 - c \phi_1) r^3 + r (2 c \phi_1 - 1)
    \end{pmatrix} \begin{pmatrix}
        \sqrt{1 - r^2} c \phi_2 \\
        \sqrt{1 - r^2} s \phi_2 \\
        r
    \end{pmatrix} \\
    &= -(1 - r^2) c \phi_2 + 2 (1 - c \phi_1) r^2 (1 - r^2) c \phi_2 + 2 r (1 - r^2) s \phi_1 s \phi_2 + 2 (1 - c \phi_1) r^4 + 2 r^2 c \phi_1 - r^2 \\
    &= \left((1 - r^2) c \phi_2 + r^2 \right) (-1 + 2 r^2) + 2 r^2 (1 - r^2) (1 - c \phi_2) c \phi_1 + 2 r (1 - r^2) s \phi_2 s \phi_1.
\end{align*}
The RHS of Eq.~\eqref{eq: uRTLGLuL} can be expanded as
\begin{align*}
    &\mathbf{u}^T_{R} \begin{pmatrix}
        \alpha_{11} & \alpha_{12} & \alpha_{13} \\
        \alpha_{21} & \alpha_{22} & \alpha_{23} \\
        \alpha_{31} & \alpha_{32} & \alpha_{33}
    \end{pmatrix} \mathbf{u}_{L} \\
    &= \begin{pmatrix}
        -\sqrt{1 - r^2} & 0 & r
    \end{pmatrix} \begin{pmatrix}
        \alpha_{11} & \alpha_{12} & \alpha_{13} \\
        \alpha_{21} & \alpha_{22} & \alpha_{23} \\
        \alpha_{31} & \alpha_{32} & \alpha_{33}
    \end{pmatrix} \begin{pmatrix}
        \sqrt{1 - r^2} \\
        0 \\
        r
    \end{pmatrix} \\
    &= \begin{pmatrix}
        -\sqrt{1 - r^2} \alpha_{11} + r \alpha_{31} & -\sqrt{1 - r^2} \alpha_{12} + r \alpha_{32} & -\sqrt{1 - r^2} \alpha_{13} + r \alpha_{33}
    \end{pmatrix} \begin{pmatrix}
        \sqrt{1 - r^2} \\
        0 \\
        r
    \end{pmatrix} \\
    &= -(1 - r^2) \alpha_{11} + r \sqrt{1 - r^2} (\alpha_{31} - \alpha_{13}) + r^2 \alpha_{33} = - \alpha_{11} + r^2 (\alpha_{11} + \alpha_{33}) + r \sqrt{1 - r^2} (\alpha_{31} - \alpha_{13}).
\end{align*}
Substituting the obtained expressions in Eq.~\eqref{eq: uRTLGLuL} and using the expression for $(1 - r^2) c \phi_2 + r^2$ from Eq.~\eqref{eq: cos_phi_2_LGL_path_expression},
\begin{align*}
    2 r^2 (1 - r^2) (1 - c \phi_2) c \phi_1 + 2 r (1 - r^2) s \phi_2 s \phi_1 &= - \alpha_{11} + r^2 (\alpha_{11} + \alpha_{33}) + r \sqrt{1 - r^2} (\alpha_{31} - \alpha_{13}) \\
    & \quad\, -  \left((1 - r^2) c \phi_2 + r^2 \right) (-1 + 2 r^2) \\
    &= - \alpha_{11} + r^2 (\alpha_{11} + \alpha_{33}) + r \sqrt{1 - r^2} (\alpha_{31} - \alpha_{13}) \\
    & \quad\, - \left(\alpha_{11} + r \sqrt{1 - r^2} (\alpha_{13} + \alpha_{31}) + r^2 (\alpha_{33} - \alpha_{11}) \right) (-1 + 2 r^2) \\
    &= 2 (\alpha_{33} - \alpha_{11}) r^2 (1 - r^2) - 2 \alpha_{13} r^3 \sqrt{1 - r^2} + 2 \alpha_{31} r (1 - r^2)^{\frac{3}{2}}. \\
    \implies r (1 - c \phi_2) c \phi_1 + s \phi_2 s \phi_1 &= (\alpha_{33} - \alpha_{11}) r - \alpha_{13} r^2 (1 - r^2)^{-\frac{1}{2}} + \alpha_{31} \sqrt{1 - r^2}.
\end{align*}
It should be noted that if $\phi_2 \neq 0,$ then the above equation can be used to solve for $\phi_1.$ Diving both sides by $\sqrt{r^2 (1 - c \phi_2)^2 + s^2 \phi_2}$ and defining $c \beta := \frac{r (1 - c \phi_2)}{\sqrt{r^2 (1 - c \phi_2)^2 + s^2 \phi_2}}, s \beta := \frac{s \phi_2}{\sqrt{r^2 (1 - c \phi_2)^2 + s^2 \phi_2}},$ the above equation can be rewritten as
\begin{align}
    c (\phi_1 - \beta) &= \frac{(\alpha_{33} - \alpha_{11}) r - \alpha_{13} r^2 (1 - r^2)^{-\frac{1}{2}} + \alpha_{31} \sqrt{1 - r^2}}{\sqrt{r^2 (1 - c \phi_2)^2 + s^2 \phi_2}}. \\
\begin{split}
    \therefore \phi_1 &= \left(\cos^{-1} \left(\frac{(\alpha_{33} - \alpha_{11}) r - \alpha_{13} r^2 (1 - r^2)^{-\frac{1}{2}} + \alpha_{31} \sqrt{1 - r^2}}{\sqrt{r^2 (1 - c \phi_2)^2 + s^2 \phi_2}} \right) + \beta \right) \% 2 \pi \\
    &= \left(\cos^{-1} \left(\frac{(\alpha_{33} - \alpha_{11}) r - \alpha_{13} r^2 (1 - r^2)^{-\frac{1}{2}} + \alpha_{31} \sqrt{1 - r^2}}{\sqrt{r^2 (1 - c \phi_2)^2 + s^2 \phi_2}} \right) + \text{atan2} \left(s \phi_2, r (1 - c \phi_2) \right) \right) \% 2 \pi.
\end{split}
\end{align}
\textbf{Important remark:} From the above equation, for each solution of $\phi_2,$ at most two solutions can be obtained for $\phi_1,$ since $\cos^{-1}$ can be replaced with $2 \pi - \cos^{-1}$ to obtain another solution. The same argument applies for the rest of the report whenever we encounter a solution using $\cos^{-1}.$ Furthermore, for all angles, a modulus function is utilized to ensure that it is between $0$ to $2 \pi.$

Finally, consider pre-multiplying Eq.~\eqref{eq: LGL_general_path_solve} by $\mathbf{u}_{L}^T$ and post-multiplying by $\mathbf{u}_{R},$ which yields
\begin{align} \label{eq: uLTLGLuR}
    \mathbf{u}^T_{L} R_L (r, \phi_1) R_G (\phi_2) R_L (r, \phi_3) \mathbf{u}_{R} = \mathbf{u}^T_{L} R_G (\phi_2) R_L (r, \phi_3) \mathbf{u}_{R} = \mathbf{u}^T_{L} \begin{pmatrix}
        \alpha_{11} & \alpha_{12} & \alpha_{13} \\
        \alpha_{21} & \alpha_{22} & \alpha_{23} \\
        \alpha_{31} & \alpha_{32} & \alpha_{33}
    \end{pmatrix} \mathbf{u}_{R}.
\end{align}
The LHS of the above equation can be expanded as
\begin{align*}
    &\mathbf{u}^T_{L} R_G (\phi_2) R_L (r, \phi_3) \mathbf{u}_{R} \\
    &= \begin{pmatrix}
        \sqrt{1 - r^2} & 0 & r
    \end{pmatrix} \begin{pmatrix}
        c \phi_2 & - s \phi_2 & 0 \\
        s \phi_2 & c \phi_2 & 0 \\
        0 & 0 & 1
    \end{pmatrix} R_L (r, \phi_3) \mathbf{u}_{R} \\
    &= \begin{pmatrix}
        \sqrt{1 - r^2} c \phi_2 & -\sqrt{1 - r^2} s \phi_2 & r
    \end{pmatrix} \begin{pmatrix}
        1 - (1 - c \phi_3) r^2 & - r s \phi_3 & (1 - c \phi_3) r \sqrt{1 - r^2} \\
        r s \phi_3 & c \phi_3 & - s \phi_3 \sqrt{1 - r^2} \\
        (1 - c \phi_3) r \sqrt{1 - r^2} & s \phi_3 \sqrt{1 - r^2} & c \phi_3 + (1 - c \phi_3) r^2
    \end{pmatrix} \begin{pmatrix}
        -\sqrt{1 - r^2} \\
        0 \\
        r
    \end{pmatrix} \\
    &= \begin{pmatrix}
        \sqrt{1 - r^2} c \phi_2 & -\sqrt{1 - r^2} s \phi_2 & r
    \end{pmatrix} \begin{pmatrix}
        -\sqrt{1 - r^2} + 2 (1 - c \phi_3) r^2 \sqrt{1 - r^2} \\
        -2 r \sqrt{1 - r^2} s \phi_3 \\
        2 (1 - c \phi_3) r^3 + r (2 c \phi_3 - 1)
    \end{pmatrix} \\
    &= - (1 - r^2) c \phi_2 + 2 (1 - c \phi_3) r^2 (1 - r^2) c \phi_2 + 2 r (1 - r^2) s \phi_3 s \phi_2 + 2 (1 - c \phi_3) r^4 + 2 r^2 c \phi_3 - r^2 \\
    &= \left((1 - r^2) c \phi_2 + r^2 \right) (-1 + 2 r^2) + 2 r^2 (1 - r^2) (1 - c \phi_2) c \phi_3 + 2 r (1 - r^2) s \phi_2 s \phi_3.
\end{align*}
The RHS of Eq.~\eqref{eq: uLTLGLuR} can be expanded as
\begin{align*}
    &\mathbf{u}^T_{L} \begin{pmatrix}
        \alpha_{11} & \alpha_{12} & \alpha_{13} \\
        \alpha_{21} & \alpha_{22} & \alpha_{23} \\
        \alpha_{31} & \alpha_{32} & \alpha_{33}
    \end{pmatrix} \mathbf{u}_{R} \\
    &= \begin{pmatrix}
        \sqrt{1 - r^2} & 0 & r
    \end{pmatrix} \begin{pmatrix}
        \alpha_{11} & \alpha_{12} & \alpha_{13} \\
        \alpha_{21} & \alpha_{22} & \alpha_{23} \\
        \alpha_{31} & \alpha_{32} & \alpha_{33}
    \end{pmatrix} \begin{pmatrix}
        -\sqrt{1 - r^2} \\
        0 \\
        r
    \end{pmatrix} \\
    &= \begin{pmatrix}
        \sqrt{1 - r^2} \alpha_{11} + r \alpha_{31} & \sqrt{1 - r^2} \alpha_{12} + r \alpha_{32} & \sqrt{1 - r^2} \alpha_{13} + r \alpha_{33}
    \end{pmatrix} \begin{pmatrix}
        -\sqrt{1 - r^2} \\
        0 \\
        r
    \end{pmatrix} \\
    &= - (1 - r^2) \alpha_{11} + r \sqrt{1 - r^2} (\alpha_{13} - \alpha_{31}) + r^2 \alpha_{33} = - \alpha_{11} + r^2 (\alpha_{11} + \alpha_{33}) + r \sqrt{1 - r^2} (\alpha_{13} - \alpha_{31}).
\end{align*}
Substituting the obtained expressions in Eq.~\eqref{eq: uLTLGLuR} and using the expression for $(1 - r^2) c \phi_2 + r^2$ from Eq.~\eqref{eq: cos_phi_2_LGL_path_expression},
\begin{align*}
    2 r^2 (1 - r^2) (1 - c \phi_2) c \phi_3 + 2 r (1 - r^2) s \phi_2 s \phi_3 &= - \alpha_{11} + r^2 (\alpha_{11} + \alpha_{33}) + r \sqrt{1 - r^2} (\alpha_{13} - \alpha_{31}) \\
    & \quad\, - \left((1 - r^2) c \phi_2 + r^2 \right) (-1 + 2 r^2) \\
    &= - \alpha_{11} + r^2 (\alpha_{11} + \alpha_{33}) + r \sqrt{1 - r^2} (\alpha_{13} - \alpha_{31}) \\
    & \quad\, - \left(\alpha_{11} + r \sqrt{1 - r^2} (\alpha_{13} + \alpha_{31}) + r^2 (\alpha_{33} - \alpha_{11}) \right) (-1 + 2 r^2) \\
    &= 2 (\alpha_{33} - \alpha_{11}) r^2 (1 - r^2) + 2 \alpha_{13} r (1 - r^2)^{\frac{3}{2}} - 2 \alpha_{31} r^3 \sqrt{1 - r^2}. \\
    \implies r (1 - c \phi_2) c \phi_3 + s \phi_2 s \phi_3 &= (\alpha_{33} - \alpha_{11}) r + \alpha_{13} \sqrt{1 - r^2} - \alpha_{31} r^2 (1 - r^2)^{-\frac{1}{2}}.
\end{align*}
It should be noted that if $\phi_2 \neq 0,$ then the above equation can be used to solve for $\phi_3.$ Diving both sides by $\sqrt{r^2 (1 - c \phi_2)^2 + s^2 \phi_2}$ and defining $c \beta := \frac{r (1 - c \phi_2)}{\sqrt{r^2 (1 - c \phi_2)^2 + s^2 \phi_2}}, s \beta := \frac{s \phi_2}{\sqrt{r^2 (1 - c \phi_2)^2 + s^2 \phi_2}},$ the above equation can be rewritten as
\begin{align}
    c (\phi_3 - \beta) &= \frac{(\alpha_{33} - \alpha_{11}) r + \alpha_{13} \sqrt{1 - r^2} - \alpha_{31} r^2 (1 - r^2)^{-\frac{1}{2}}}{\sqrt{r^2 (1 - c \phi_2)^2 + s^2 \phi_2}}. \\
\begin{split}
    \therefore \phi_3 &= \left(\cos^{-1} \left(\frac{(\alpha_{33} - \alpha_{11}) r + \alpha_{13} \sqrt{1 - r^2} - \alpha_{31} r^2 (1 - r^2)^{-\frac{1}{2}}}{\sqrt{r^2 (1 - c \phi_2)^2 + s^2 \phi_2}} \right) + \beta \right) \% 2 \pi \\
    &= \left(\cos^{-1} \left(\frac{(\alpha_{33} - \alpha_{11}) r + \alpha_{13} \sqrt{1 - r^2} - \alpha_{31} r^2 (1 - r^2)^{-\frac{1}{2}}}{\sqrt{r^2 (1 - c \phi_2)^2 + s^2 \phi_2}} \right) + \tan^{-1} \left(\frac{s \phi_2}{r (1 - c \phi_2)} \right) \right) \% 2 \pi.
\end{split}
\end{align}
From the above equation, for each solution of $\phi_2,$ at most two solutions can be obtained for $\phi_3$.

\textbf{Remark:} In the case that $\phi_2 = 0,$ the $LGL$ path reduces to an $L$ segment. Hence, in this case, $\phi_3$ can be set to zero without loss of generality, and $\phi_1$ can be solved using the first and second entries in the second row in the net rotation matrix (refer to the expression for $R_L$). That is,
\begin{align} \label{eq: expression_phi1_when_phi2_zero}
    \phi_1 = \left(\text{atan2} \left(\alpha_{21}, r \alpha_{22} \right) \right) \% 2 \pi.
\end{align}

\subsection{RGR path}

The equation to be solved is given by
\begin{align} \label{appeq: RGR_general_path_solve}
    R_R (r, \phi_1) R_G (\phi_2) R_R (r, \phi_3) = \begin{pmatrix}
        \alpha_{11} & \alpha_{12} & \alpha_{13} \\
        \alpha_{21} & \alpha_{22} & \alpha_{23} \\
        \alpha_{31} & \alpha_{32} & \alpha_{33}
    \end{pmatrix}.
\end{align}

Pre-multiplying Eq.~\eqref{appeq: RGR_general_path_solve} by $\mathbf{u}^T_R$ and post-multiplying by $\mathbf{u}_R,$
\begin{align*}
    \mathbf{u}^T_R R_R (r, \phi_1) R_G (\phi_2) R_R (r, \phi_3) \mathbf{u}_R = \mathbf{u}^T_R R_G (\phi_2) \mathbf{u}_R &= \mathbf{u}^T_R \begin{pmatrix}
        \alpha_{11} & \alpha_{12} & \alpha_{13} \\
        \alpha_{21} & \alpha_{22} & \alpha_{23} \\
        \alpha_{31} & \alpha_{32} & \alpha_{33}
    \end{pmatrix} \mathbf{u}_R. \\
    \therefore \begin{pmatrix}
        -\sqrt{1 - r^2} & 0 & r
    \end{pmatrix} \begin{pmatrix}
        c \phi_2 & - s \phi_2 & 0 \\
        s \phi_2 & c \phi_2 & 0 \\
        0 & 0 & 1
    \end{pmatrix} \begin{pmatrix}
        -\sqrt{1 - r^2} \\
        0 \\
        r
    \end{pmatrix} &= \begin{pmatrix}
        -\sqrt{1 - r^2} & 0 & r
    \end{pmatrix} \begin{pmatrix}
        \alpha_{11} & \alpha_{12} & \alpha_{13} \\
        \alpha_{21} & \alpha_{22} & \alpha_{23} \\
        \alpha_{31} & \alpha_{32} & \alpha_{33}
    \end{pmatrix} \begin{pmatrix}
        -\sqrt{1 - r^2} \\
        0 \\
        r
    \end{pmatrix}. 
    % \\
    % \implies \begin{pmatrix}
    %     -\sqrt{1 - r^2} c \phi_2 & \sqrt{1 - r^2} s \phi_2 & r
    % \end{pmatrix} \begin{pmatrix}
    %     -\sqrt{1 - r^2} \\
    %     0 \\
    %     r
    % \end{pmatrix} &= \begin{pmatrix}
    %     -\sqrt{1 - r^2} \alpha_{11} + r \alpha_{31} & -\sqrt{1 - r^2} \alpha_{12} + r \alpha_{32} & -\sqrt{1 - r^2} \alpha_{13} + r \alpha_{33} 
    % \end{pmatrix}
\end{align*}
Simplifying the above equation,
\begin{align} \label{appeq: cos_phi_2_RGR_path_expression}
\begin{split}
    (1 - r^2) c \phi_2 + r^2 &= (1 - r^2) \alpha_{11} - r \sqrt{1 - r^2} (\alpha_{13} + \alpha_{31}) + r^2 \alpha_{33}. \\
    % \implies c \phi_2 &= \frac{(1 - r^2) \alpha_{11} - r \sqrt{1 - r^2} (\alpha_{13} + \alpha_{31}) + r^2 \alpha_{33} - r^2}{1 - r^2}, \\
    \implies c \phi_2 &= \frac{\alpha_{11} - r \sqrt{1 - r^2} (\alpha_{13} + \alpha_{31}) + r^2 (\alpha_{33} - \alpha_{11} - 1)}{1 - r^2},
\end{split}
\end{align}
which yields at most two solutions for $\phi_2 \in [0, 2 \pi)$ if the absolute value of the RHS is less than or equal to 1.

Consider pre-multiplying Eq.~\eqref{appeq: RGR_general_path_solve} by $\mathbf{u}_L^T$ and post-multiplying by $\mathbf{u}_R$, which yields
\begin{align} \label{appeq: uLTRGRuR}
    \mathbf{u}^T_L R_R (r, \phi_1) R_G (\phi_2) R_R (r, \phi_3) \mathbf{u}_R = \mathbf{u}^T_L R_R (r, \phi_1) R_G (\phi_2) \mathbf{u}_R = \mathbf{u}^T_L \begin{pmatrix}
        \alpha_{11} & \alpha_{12} & \alpha_{13} \\
        \alpha_{21} & \alpha_{22} & \alpha_{23} \\
        \alpha_{31} & \alpha_{32} & \alpha_{33}
    \end{pmatrix} \mathbf{u}_R.
\end{align}
The LHS of the above equation can be expanded as
\begin{align*}
    &\mathbf{u}^T_L R_R (r, \phi_1) R_G (\phi_2) \mathbf{u}_R \\
    &= \mathbf{u}^T_L R_R (r, \phi_1) \begin{pmatrix}
        c \phi_2 & - s \phi_2 & 0 \\
        s \phi_2 & c \phi_2 & 0 \\
        0 & 0 & 1
    \end{pmatrix} \begin{pmatrix}
        -\sqrt{1 - r^2} \\
        0 \\
        r
    \end{pmatrix} \\
    &= \begin{pmatrix}
        \sqrt{1 - r^2} & 0 & r
    \end{pmatrix} \begin{pmatrix}
        1 - (1 - c \phi_1) r^2 & - r s \phi_1 & -(1 - c \phi_1) r \sqrt{1 - r^2} \\
        r s \phi_1 & c \phi_1 & s \phi_1 \sqrt{1 - r^2} \\
        -(1 - c \phi_1) r \sqrt{1 - r^2} & -s \phi_1 \sqrt{1 - r^2} & c \phi_1 + (1 - c \phi_1) r^2
    \end{pmatrix} \begin{pmatrix}
        -\sqrt{1 - r^2} c \phi_2 \\
        -\sqrt{1 - r^2} s \phi_2 \\
        r
    \end{pmatrix} \\
    &= \begin{pmatrix}
        \sqrt{1 - r^2} - 2 (1 - c \phi_1) r^2 \sqrt{1 - r^2} & -2 r \sqrt{1 - r^2} s \phi_1 & 2 (1 - c \phi_1) r^3 + r (2 c \phi_1 - 1)
    \end{pmatrix} \begin{pmatrix}
        -\sqrt{1 - r^2} c \phi_2 \\
        -\sqrt{1 - r^2} s \phi_2 \\
        r
    \end{pmatrix} \\
    &= - (1 - r^2) c \phi_2 + 2 (1 - c \phi_1) r^2 (1 - r^2) c \phi_2 + 2 r (1 - r^2) s \phi_1 s \phi_2 + 2 (1 - c \phi_1) r^4 + 2 r^2 c \phi_1 - r^2 \\
    &= \left((1 - r^2) c \phi_2 + r^2 \right) (-1 + 2 r^2) + 2 r^2 (1 - r^2) (1 - c \phi_2) c \phi_1 + 2 r (1 - r^2) s \phi_2 s \phi_1.
\end{align*}
The RHS of Eq.~\eqref{appeq: uLTRGRuR} can be expanded as
\begin{align*}
    &\mathbf{u}^T_L \begin{pmatrix}
        \alpha_{11} & \alpha_{12} & \alpha_{13} \\
        \alpha_{21} & \alpha_{22} & \alpha_{23} \\
        \alpha_{31} & \alpha_{32} & \alpha_{33}
    \end{pmatrix} \mathbf{u}_R \\
    &= \begin{pmatrix}
        \sqrt{1 - r^2} & 0 & r
    \end{pmatrix} \begin{pmatrix}
        \alpha_{11} & \alpha_{12} & \alpha_{13} \\
        \alpha_{21} & \alpha_{22} & \alpha_{23} \\
        \alpha_{31} & \alpha_{32} & \alpha_{33}
    \end{pmatrix} \begin{pmatrix}
        -\sqrt{1 - r^2} \\
        0 \\
        r
    \end{pmatrix} \\
    &= \begin{pmatrix}
        \sqrt{1 - r^2} \alpha_{11} + r \alpha_{31} & \sqrt{1 - r^2} \alpha_{12} + r \alpha_{32} & \sqrt{1 - r^2} \alpha_{13} + r \alpha_{33}
    \end{pmatrix} \begin{pmatrix}
        -\sqrt{1 - r^2} \\
        0 \\
        r
    \end{pmatrix} \\
    &= - (1 - r^2) \alpha_{11} + r \sqrt{1 - r^2} (\alpha_{13} - \alpha_{31}) + r^2 \alpha_{33} = - \alpha_{11} + r^2 (\alpha_{11} + \alpha_{33}) + r \sqrt{1 - r^2} (\alpha_{13} - \alpha_{31}).
\end{align*}

Substituting the obtained expressions in Eq.~\eqref{appeq: uLTRGRuR} and using the expression for $(1 - r^2) c \phi_2 + r^2$ from Eq.~\eqref{appeq: cos_phi_2_RGR_path_expression},
\begin{align*}
    2 r^2 (1 - r^2) (1 - c \phi_2) c \phi_1 + 2 r (1 - r^2) s \phi_2 s \phi_1 &= - \alpha_{11} + r^2 (\alpha_{11} + \alpha_{33}) + r \sqrt{1 - r^2} (\alpha_{13} - \alpha_{31}) \\
    &\quad\, - \left((1 - r^2) c \phi_2 + r^2 \right) (-1 + 2 r^2) \\
    &= - \alpha_{11} + r^2 (\alpha_{11} + \alpha_{33}) + r \sqrt{1 - r^2} (\alpha_{13} - \alpha_{31}) \\
    &\quad\, - \left(\alpha_{11} + r^2 (\alpha_{33} - \alpha_{11}) - r \sqrt{1 - r^2} (\alpha_{13} + \alpha_{31}) \right) (-1 + 2 r^2) \\
    &= 2 (\alpha_{33} - \alpha_{11}) r^2 (1 - r^2) + 2 \alpha_{13} r^3 \sqrt{1 - r^2} - 2 \alpha_{31} r (1 - r^2)^{\frac{3}{2}}. \\
    \implies r (1 - c \phi_2) c \phi_1 + s \phi_2 s \phi_1 &= (\alpha_{33} - \alpha_{11}) r + \alpha_{13} r^2 (1 - r^2)^{-\frac{1}{2}} - \alpha_{31} (1 - r^2)^{\frac{1}{2}}.
\end{align*}
It should be noted that if $\phi_2 \neq 0,$ then the above equation can be used to solve for $\phi_1.$ Diving both sides by $\sqrt{r^2 (1 - c \phi_2)^2 + s^2 \phi_2}$ and defining $c \beta := \frac{r (1 - c \phi_2)}{\sqrt{r^2 (1 - c \phi_2)^2 + s^2 \phi_2}}, s \beta := \frac{s \phi_2}{\sqrt{r^2 (1 - c \phi_2)^2 + s^2 \phi_2}},$ the above equation can be rewritten as
\begin{align}
    c (\phi_1 - \beta) &= \frac{(\alpha_{33} - \alpha_{11}) r + \alpha_{13} r^2 (1 - r^2)^{-\frac{1}{2}} - \alpha_{31} (1 - r^2)^{\frac{1}{2}}}{\sqrt{r^2 (1 - c \phi_2)^2 + s^2 \phi_2}}. \\
\begin{split}
    \therefore \phi_1 &= \left(\cos^{-1} \left(\frac{(\alpha_{33} - \alpha_{11}) r + \alpha_{13} r^2 (1 - r^2)^{-\frac{1}{2}} - \alpha_{31} (1 - r^2)^{\frac{1}{2}}}{\sqrt{r^2 (1 - c \phi_2)^2 + s^2 \phi_2}} \right) + \beta \right) \% 2 \pi \\
    &= \left(\cos^{-1} \left(\frac{(\alpha_{33} - \alpha_{11}) r + \alpha_{13} r^2 (1 - r^2)^{-\frac{1}{2}} - \alpha_{31} (1 - r^2)^{\frac{1}{2}}}{\sqrt{r^2 (1 - c \phi_2)^2 + s^2 \phi_2}} \right) + \text{atan2} \left(s \phi_2, r (1 - c \phi_2) \right) \right) \% 2 \pi.
\end{split}
\end{align}
From the above equation, for each solution of $\phi_2,$ at most two solutions can be obtained for $\phi_1$.

Finally, consider pre-multiplying Eq.~\eqref{appeq: RGR_general_path_solve} by $\mathbf{u}_R^T$ and post-multiplying by $\mathbf{u}_L,$ which yields
\begin{align} \label{appeq: uRTRGRuL}
    \mathbf{u}^T_R R_R (r, \phi_1) R_G (\phi_2) R_R (r, \phi_3) \mathbf{u}_L = \mathbf{u}^T_R R_G (\phi_2) R_R (r, \phi_3) \mathbf{u}_L = \mathbf{u}^T_R \begin{pmatrix}
        \alpha_{11} & \alpha_{12} & \alpha_{13} \\
        \alpha_{21} & \alpha_{22} & \alpha_{23} \\
        \alpha_{31} & \alpha_{32} & \alpha_{33}
    \end{pmatrix} \mathbf{u}_L.
\end{align}
The LHS of the above equation can be expanded as
\begin{align*}
    &\mathbf{u}^T_R R_G (\phi_2) R_R (r, \phi_3) \mathbf{u}_L \\
    &= \begin{pmatrix}
        -\sqrt{1 - r^2} & 0 & r
    \end{pmatrix} \begin{pmatrix}
        c \phi_2 & - s \phi_2 & 0 \\
        s \phi_2 & c \phi_2 & 0 \\
        0 & 0 & 1
    \end{pmatrix} R_R (r, \phi_3) \mathbf{u}_L \\
    &= \begin{pmatrix}
        -\sqrt{1 - r^2} c \phi_2 & \sqrt{1 - r^2} s \phi_2 & r
    \end{pmatrix} \begin{pmatrix}
        1 - (1 - c \phi_3) r^2 & - r s \phi_3 & -(1 - c \phi_3) r \sqrt{1 - r^2} \\
        r s \phi_3 & c \phi_3 & s \phi_3 \sqrt{1 - r^2} \\
        -(1 - c \phi_3) r \sqrt{1 - r^2} & -s \phi_3 \sqrt{1 - r^2} & c \phi_3 + (1 - c \phi_3) r^2
    \end{pmatrix} \begin{pmatrix}
        \sqrt{1 - r^2} \\
        0 \\
        r
    \end{pmatrix} \\
    &= \begin{pmatrix}
        -\sqrt{1 - r^2} c \phi_2 & \sqrt{1 - r^2} s \phi_2 & r
    \end{pmatrix} \begin{pmatrix}
        \sqrt{1 - r^2} - 2 (1 - c \phi_3) r^2 \sqrt{1 - r^2} \\
        2 r \sqrt{1 - r^2} s \phi_3 \\
        2 (1 - c \phi_3) r^3 + r (2 c \phi_3 - 1)
    \end{pmatrix} \\
    &= - (1 - r^2) c \phi_2 + 2 (1 - c \phi_3) r^2 (1 - r^2) c \phi_2 + 2 r (1 - r^2) s \phi_3 s \phi_2 + 2 (1 - c \phi_3) r^4 + 2 r^2 c \phi_3 - r^2 \\
    &= \left((1 - r^2) c \phi_2 + r^2 \right) (-1 + 2 r^2) + 2 r^2 (1 - r^2) (1 - c \phi_2) c \phi_3 + 2 r (1 - r^2) s \phi_2 s \phi_3.
\end{align*}
The RHS of Eq.~\eqref{appeq: uRTRGRuL} can be expanded as
\begin{align*}
    &\mathbf{u}^T_R \begin{pmatrix}
        \alpha_{11} & \alpha_{12} & \alpha_{13} \\
        \alpha_{21} & \alpha_{22} & \alpha_{23} \\
        \alpha_{31} & \alpha_{32} & \alpha_{33}
    \end{pmatrix} \mathbf{u}_L \\
    &= \begin{pmatrix}
        -\sqrt{1 - r^2} & 0 & r
    \end{pmatrix} \begin{pmatrix}
        \alpha_{11} & \alpha_{12} & \alpha_{13} \\
        \alpha_{21} & \alpha_{22} & \alpha_{23} \\
        \alpha_{31} & \alpha_{32} & \alpha_{33}
    \end{pmatrix} \begin{pmatrix}
        \sqrt{1 - r^2} \\
        0 \\
        r
    \end{pmatrix} \\
    &= \begin{pmatrix}
        -\sqrt{1 - r^2} \alpha_{11} + r \alpha_{31} & -\sqrt{1 - r^2} \alpha_{12} + r \alpha_{32} & -\sqrt{1 - r^2} \alpha_{13} + r \alpha_{33}
    \end{pmatrix} \begin{pmatrix}
        \sqrt{1 - r^2} \\
        0 \\
        r
    \end{pmatrix} \\
    &= - (1 - r^2) \alpha_{11} + r \sqrt{1 - r^2} (\alpha_{31} - \alpha_{13}) + r^2 \alpha_{33} = - \alpha_{11} + r^2 (\alpha_{11} + \alpha_{33}) + r \sqrt{1 - r^2} (\alpha_{31} - \alpha_{13}).
\end{align*}
Substituting the obtained expressions in Eq.~\eqref{appeq: uRTRGRuL} and using the expression for $(1 - r^2) c \phi_2 + r^2$ from Eq.~\eqref{appeq: cos_phi_2_RGR_path_expression},
\begin{align*}
    2 r^2 (1 - r^2) (1 - c \phi_2) c \phi_3 + 2 r (1 - r^2) s \phi_2 s \phi_3 &= - \alpha_{11} + r^2 (\alpha_{11} + \alpha_{33}) + r \sqrt{1 - r^2} (\alpha_{31} - \alpha_{13}) \\
    &\quad\, - \left((1 - r^2) c \phi_2 + r^2 \right) (-1 + 2 r^2) \\
    &= - \alpha_{11} + r^2 (\alpha_{11} + \alpha_{33}) + r \sqrt{1 - r^2} (\alpha_{31} - \alpha_{13}) \\
    &\quad\, - \left(\alpha_{11} + r^2 (\alpha_{33} - \alpha_{11}) - r \sqrt{1 - r^2} (\alpha_{13} + \alpha_{31}) \right) (-1 + 2 r^2) \\
    &= 2 (\alpha_{33} - \alpha_{11}) r^2 (1 - r^2) - 2 \alpha_{13} r (1 - r^2)^{\frac{3}{2}} + 2 \alpha_{31} r^3 \sqrt{1 - r^2}. \\
    \implies r (1 - c \phi_2) c \phi_3 + s \phi_2 s \phi_3 &= (\alpha_{33} - \alpha_{11}) r - \alpha_{13} (1 - r^2)^{\frac{1}{2}} + \alpha_{31} r^2 (1 - r^2)^{-\frac{1}{2}}.
\end{align*}
It should be noted that if $\phi_2 \neq 0,$ then the above equation can be used to solve for $\phi_3.$ Diving both sides by $\sqrt{r^2 (1 - c \phi_2)^2 + s^2 \phi_2}$ and defining $c \beta := \frac{r (1 - c \phi_2)}{\sqrt{r^2 (1 - c \phi_2)^2 + s^2 \phi_2}}, s \beta := \frac{s \phi_2}{\sqrt{r^2 (1 - c \phi_2)^2 + s^2 \phi_2}},$ the above equation can be rewritten as
\begin{align}
    c (\phi_3 - \beta) &= \frac{(\alpha_{33} - \alpha_{11}) r - \alpha_{13} (1 - r^2)^{\frac{1}{2}} + \alpha_{31} r^2 (1 - r^2)^{-\frac{1}{2}}}{\sqrt{r^2 (1 - c \phi_2)^2 + s^2 \phi_2}}. \\
\begin{split}
    \therefore \phi_3 &= \left(\cos^{-1} \left(\frac{(\alpha_{33} - \alpha_{11}) r - \alpha_{13} (1 - r^2)^{\frac{1}{2}} + \alpha_{31} r^2 (1 - r^2)^{-\frac{1}{2}}}{\sqrt{r^2 (1 - c \phi_2)^2 + s^2 \phi_2}} \right) + \beta \right) \% 2 \pi \\
    &= \left(\cos^{-1} \left(\frac{(\alpha_{33} - \alpha_{11}) r - \alpha_{13} (1 - r^2)^{\frac{1}{2}} + \alpha_{31} r^2 (1 - r^2)^{-\frac{1}{2}}}{\sqrt{r^2 (1 - c \phi_2)^2 + s^2 \phi_2}} \right) + \text{atan2} \left(s \phi_2, r (1 - c \phi_2) \right) \right) \% 2 \pi.
\end{split}
\end{align}
From the above equation, for each solution of $\phi_2,$ at most two solutions can be obtained for $\phi_3$.

\textbf{Remark}: The case of $\phi_2 = 0$ can be handled similar to the $LGL$ path, wherein $\phi_3$ can be set to zero and the expression for $\phi_1$ obtained is the same as that given in Eq.~\eqref{eq: expression_phi1_when_phi2_zero}.

\subsection{LGR path}

The equation to be solved is given by
\begin{align} \label{eq: LGR_general_path_solve}
    R_L (r, \phi_1) R_G (\phi_2) R_R (r, \phi_3) = \begin{pmatrix}
        \alpha_{11} & \alpha_{12} & \alpha_{13} \\
        \alpha_{21} & \alpha_{22} & \alpha_{23} \\
        \alpha_{31} & \alpha_{32} & \alpha_{33}
    \end{pmatrix}.
\end{align}
Pre-multiplying Eq.~\eqref{eq: LGR_general_path_solve} by $\mathbf{u}^T_{L}$ and post-multiplying by $\mathbf{u}_{R},$
\begin{align*}
    \mathbf{u}^T_{L} R_L (r, \phi_1) R_G (\phi_2) R_R (r, \phi_3) \mathbf{u}_{R} = \mathbf{u}^T_{L} R_G (\phi_2) \mathbf{u}_{R} &= \mathbf{u}^T_{L} \begin{pmatrix}
        \alpha_{11} & \alpha_{12} & \alpha_{13} \\
        \alpha_{21} & \alpha_{22} & \alpha_{23} \\
        \alpha_{31} & \alpha_{32} & \alpha_{33}
    \end{pmatrix} \mathbf{u}_{R}. \\
    \therefore \begin{pmatrix}
        \sqrt{1 - r^2} & 0 & r
    \end{pmatrix} \begin{pmatrix}
        c \phi_2 & - s \phi_2 & 0 \\
        s \phi_2 & c \phi_2 & 0 \\
        0 & 0 & 1
    \end{pmatrix} \begin{pmatrix}
        -\sqrt{1 - r^2} \\
        0 \\
        r
    \end{pmatrix} &= \begin{pmatrix}
        \sqrt{1 - r^2} & 0 & r
    \end{pmatrix} \begin{pmatrix}
        \alpha_{11} & \alpha_{12} & \alpha_{13} \\
        \alpha_{21} & \alpha_{22} & \alpha_{23} \\
        \alpha_{31} & \alpha_{32} & \alpha_{33}
    \end{pmatrix} \begin{pmatrix}
        -\sqrt{1 - r^2} \\
        0 \\
        r
    \end{pmatrix}.
\end{align*}
Simplifying the above equation,
\begin{align} \label{eq: cos_phi_2_LGR_path_expression}
\begin{split}
    &-\left(1 - r^2 \right) c \phi_2 + r^2 = -\left(1 - r^2 \right) \alpha_{11} + r \sqrt{1 - r^2} \left(\alpha_{13} - \alpha_{31} \right) + r^2 \alpha_{33}. \\
    &\implies c \phi_2 = \frac{\left(1 - r^2 \right) \alpha_{11} + r \sqrt{1 - r^2} \left(\alpha_{31} -  \alpha_{13} \right) + r^2 (1 - \alpha_{33})}{\left(1 - r^2 \right)},
\end{split}
\end{align}
which yields at most two solutions for $\phi_2 \in [0, 2 \pi)$ if the absolute value of the RHS is less than or equal to 1.

Consider pre-multiplying Eq.~\eqref{eq: LGR_general_path_solve} by $\mathbf{u}_{R}^T$ and post-multiplying by $\mathbf{u}_{R}$, which yields
\begin{align} \label{eq: uRLTLGRuRR}
    \mathbf{u}^T_{R} R_L (r, \phi_1) R_G (\phi_2) R_R (r, \phi_3) \mathbf{u}_{R} = \mathbf{u}^T_{R} R_L (r, \phi_1) R_G (\phi_2) \mathbf{u}_{R} = \mathbf{u}^T_{R} \begin{pmatrix}
        \alpha_{11} & \alpha_{12} & \alpha_{13} \\
        \alpha_{21} & \alpha_{22} & \alpha_{23} \\
        \alpha_{31} & \alpha_{32} & \alpha_{33}
    \end{pmatrix} \mathbf{u}_{R}.
\end{align}
The LHS of the above equation can be expanded as
\begin{align*}
    &\mathbf{u}^T_{R} R_L (r, \phi_1) R_G (\phi_2) \mathbf{u}_{R} \\
    &= \mathbf{u}^T_{R} R_L (r, \phi_1) \begin{pmatrix}
        c \phi_2 & - s \phi_2 & 0 \\
        s \phi_2 & c \phi_2 & 0 \\
        0 & 0 & 1
    \end{pmatrix} \begin{pmatrix}
        -\sqrt{1 - r^2} \\
        0 \\
        r
    \end{pmatrix} \\
    &= \begin{pmatrix}
        -\sqrt{1 - r^2} & 0 & r
    \end{pmatrix} \begin{pmatrix}
        1 - (1 - c \phi_1) r^2 & - r s \phi_1 & (1 - c \phi_1) r \sqrt{1 - r^2} \\
        r s \phi_1 & c \phi_1 & - s \phi_1 \sqrt{1 - r^2} \\
        (1 - c \phi_1) r \sqrt{1 - r^2} & s \phi_1 \sqrt{1 - r^2} & c \phi_1 + (1 - c \phi_1) r^2
    \end{pmatrix} \begin{pmatrix}
        -\sqrt{1 - r^2} c \phi_2 \\
        -\sqrt{1 - r^2} s \phi_2 \\
        r
    \end{pmatrix} \\
    &= (1 - 2 r^2) \left(\left(1 - r^2 \right) c \phi_2 - r^2 \right) + 2 r \sqrt{1 - r^2} \left(r \sqrt{1 - r^2} c \phi_2 + r \sqrt{1 - r^2} \right) c \phi_1 - 2 r \left(1 - r^2 \right) s \phi_2 s \phi_1.
\end{align*}
The RHS of Eq.~\eqref{eq: uRLTLGRuRR} can be expanded as
\begin{align*}
    \mathbf{u}^T_{R} \begin{pmatrix}
        \alpha_{11} & \alpha_{12} & \alpha_{13} \\
        \alpha_{21} & \alpha_{22} & \alpha_{23} \\
        \alpha_{31} & \alpha_{32} & \alpha_{33}
    \end{pmatrix} \mathbf{u}_{R} &= \begin{pmatrix}
        -\sqrt{1 - r^2} & 0 & r
    \end{pmatrix} \begin{pmatrix}
        \alpha_{11} & \alpha_{12} & \alpha_{13} \\
        \alpha_{21} & \alpha_{22} & \alpha_{23} \\
        \alpha_{31} & \alpha_{32} & \alpha_{33}
    \end{pmatrix} \begin{pmatrix}
        -\sqrt{1 - r^2} \\
        0 \\
        r
    \end{pmatrix} \\
    &= \left(1 - r^2 \right) \alpha_{11} - r \sqrt{1 - r^2} \alpha_{31} - r \sqrt{1 - r^2} \alpha_{13} + r^2 \alpha_{33}.
\end{align*}
Substituting the obtained expressions in Eq.~\eqref{eq: uRLTLGRuRR} and using the expression for $\left(1 - r^2 \right) c \phi_2 - r^2$ from Eq.~\eqref{eq: cos_phi_2_LGR_path_expression},
\begin{align} \label{eq: LGR_solving_phi_1}
\begin{split}
    &2 r \sqrt{1 - r^2} \left(r \sqrt{1 - r^2} c \phi_2 + r \sqrt{1 - r^2} \right) c \phi_1 - 2 r \left(1 - r^2 \right) s \phi_2 s \phi_1 \\
    &= \left(1 - r^2 \right) \alpha_{11} - r \sqrt{1 - r^2} \alpha_{31} - r \sqrt{1 - r^2} \alpha_{13} + r^2 \alpha_{33} \\
    & \quad\, - (1 - 2 r^2) \left(\left(1 - r^2 \right) \alpha_{11} + r \sqrt{1 - r^2} \alpha_{31} - r \sqrt{1 - r^2} \alpha_{13} - r^2 \alpha_{33} \right) \\
    &= 2 r^2 \left(1 - r^2 \right) \alpha_{11} - 2 r (1 - r^2) \sqrt{1 - r^2} \alpha_{31} - 2 r^3 \sqrt{1 - r^2} \alpha_{13} + 2 r^2 (1 - r^2) \alpha_{33}. \\
    \implies& \left(r \sqrt{1 - r^2} c \phi_2 + r \sqrt{1 - r^2} \right) c \phi_1 - \sqrt{1 - r^2} s \phi_2 s \phi_1 \\
    &= r \sqrt{1 - r^2} \alpha_{11} - \left(1 - r^2 \right) \alpha_{31} - r^2 \alpha_{13} + r \sqrt{1 - r^2} \alpha_{33}.
\end{split}
\end{align}
It should be noted that if $\phi_2 = \pi$, the coefficient of $s \phi_1$ and $c \phi_1$ reduces to $0.$ This special case will be treated separately. 

When $\phi_2 \neq \pi,$ dividing both sides of Eq.~\eqref{eq: LGR_solving_phi_1} by $\sqrt{\left(r \sqrt{1 - r^2} c \phi_2 + r \sqrt{1 - r^2} \right)^2 + (1 - r^2) s^2 \phi_2}$ and defining
\begin{align*}
    c \beta &:= \frac{r \sqrt{1 - r^2} c \phi_2 + r \sqrt{1 - r^2}}{\sqrt{\left(r \sqrt{1 - r^2} c \phi_2 + r \sqrt{1 - r^2} \right)^2 + (1 - r^2) s^2 \phi_2}}, \quad
    s \beta := \frac{\sqrt{1 - r^2} s \phi_2}{\sqrt{\left(r \sqrt{1 - r^2} c \phi_2 + r \sqrt{1 - r^2} \right)^2 + (1 - r^2) s^2 \phi_2}},
\end{align*}
Eq.~\eqref{eq: LGR_solving_phi_1} can be rewritten as
\begin{align}
    c (\phi_1 + \beta) &= \frac{r \sqrt{1 - r^2} \alpha_{11} - \left(1 - r^2 \right) \alpha_{31} - r^2 \alpha_{13} + r \sqrt{1 - r^2} \alpha_{33}}{\sqrt{\left(r \sqrt{1 - r^2} c \phi_2 + r \sqrt{1 - r^2} \right)^2 + (1 - r^2) s^2 \phi_2}}. \\
\begin{split}
    \therefore \phi_1 &= \left(\cos^{-1} {\left(\frac{r \sqrt{1 - r^2} \alpha_{11} - \left(1 - r^2 \right) \alpha_{31} - r^2 \alpha_{13} + r \sqrt{1 - r^2} \alpha_{33}}{\sqrt{\left(r \sqrt{1 - r^2} c \phi_2 + r \sqrt{1 - r^2} \right)^2 + (1 - r^2) s^2 \phi_2}} \right)} - \beta \right) \% 2 \pi,
\end{split}
\end{align}
where $\beta = \text{atan2} \left(s \phi_2, r \left(c \phi_2 + 1 \right) \right)$. From the above equation, for each solution of $\phi_2,$ at most two solutions can be obtained for $\phi_1$.

Consider pre-multiplying Eq.~\eqref{eq: LGR_general_path_solve} by $\mathbf{u}_{L}^T$ and post-multiplying by $\mathbf{u}_{L}$, which yields
\begin{align} \label{eq: uLLTLGRuLR}
    \mathbf{u}^T_{L} R_L (r, \phi_1) R_G (\phi_2) R_R (r, \phi_3) \mathbf{u}_{L} = \mathbf{u}^T_{L} R_G (\phi_2) R_R (r, \phi_3) \mathbf{u}_{L} = \mathbf{u}^T_{L} \begin{pmatrix}
        \alpha_{11} & \alpha_{12} & \alpha_{13} \\
        \alpha_{21} & \alpha_{22} & \alpha_{23} \\
        \alpha_{31} & \alpha_{32} & \alpha_{33}
    \end{pmatrix} \mathbf{u}_{L}.
\end{align}
The LHS of the above equation can be expanded as
\begin{align*}
    &\mathbf{u}^T_{L} R_G (\phi_2) R_R (r, \phi_3) \mathbf{u}_{L} \\
    &= \begin{pmatrix}
        \sqrt{1 - r^2} & 0 & r
    \end{pmatrix} \begin{pmatrix}
        c \phi_2 & - s \phi_2 & 0 \\
        s \phi_2 & c \phi_2 & 0 \\
        0 & 0 & 1
    \end{pmatrix} R_R (r, \phi_3) \mathbf{u}_{L} \\
    &= \begin{pmatrix}
        \sqrt{1 - r^2} c \phi_2 & -\sqrt{1 - r^2} s \phi_2 & r
    \end{pmatrix} \begin{pmatrix}
        1 - (1 - c \phi_3) r^2 & - r s \phi_3 & -(1 - c \phi_3) r \sqrt{1 - r^2} \\
        r s \phi_3 & c \phi_3 & s \phi_3 \sqrt{1 - r^2} \\
        -(1 - c \phi_3) r \sqrt{1 - r^2} & -s \phi_3 \sqrt{1 - r^2} & c \phi_3 + (1 - c \phi_3) r^2
    \end{pmatrix} \begin{pmatrix}
        \sqrt{1 - r^2} \\
        0 \\
        r
    \end{pmatrix} \\
    &= (1 - 2 r^2) \left(\left(1 - r^2 \right) c \phi_2 - r^2 \right) + 2 r \sqrt{1 - r^2} \left(r \sqrt{1 - r^2} c \phi_2 + r \sqrt{1 - r^2} \right) c \phi_3 - 2 r \left(1 - r^2 \right) s \phi_2 s \phi_3.
\end{align*}
The RHS of Eq.~\eqref{eq: uLLTLGRuLR} can be expanded as
\begin{align*}
    \mathbf{u}^T_{L} \begin{pmatrix}
        \alpha_{11} & \alpha_{12} & \alpha_{13} \\
        \alpha_{21} & \alpha_{22} & \alpha_{23} \\
        \alpha_{31} & \alpha_{32} & \alpha_{33}
    \end{pmatrix} \mathbf{u}_{L} &= \begin{pmatrix}
        \sqrt{1 - r^2} & 0 & r
    \end{pmatrix} \begin{pmatrix}
        \alpha_{11} & \alpha_{12} & \alpha_{13} \\
        \alpha_{21} & \alpha_{22} & \alpha_{23} \\
        \alpha_{31} & \alpha_{32} & \alpha_{33}
    \end{pmatrix} \begin{pmatrix}
        \sqrt{1 - r^2} \\
        0 \\
        r
    \end{pmatrix} \\
    &= \left(1 - r^2 \right) \alpha_{11} + r \sqrt{1 - r^2} \alpha_{31} + r \sqrt{1 - r^2} \alpha_{13} + r^2 \alpha_{33}.
\end{align*}
Substituting the obtained expressions in Eq.~\eqref{eq: uLLTLGRuLR} and using the expression for $\left(1 - r^2 \right) c \phi_2 - r^2$ from Eq.~\eqref{eq: cos_phi_2_LGR_path_expression},
\begin{align} \label{eq: LGR_solving_phi_3}
\begin{split}
    &2 r \sqrt{1 - r^2} \left(r \sqrt{1 - r^2} c \phi_2 + r \sqrt{1 - r^2} \right) c \phi_3 - 2 r \left(1 - r^2 \right) s \phi_2 s \phi_3 \\
    &= \left(1 - r^2 \right) \alpha_{11} + r \sqrt{1 - r^2} \alpha_{31} + r \sqrt{1 - r^2} \alpha_{13} + r^2 \alpha_{33} \\
    & \quad\, - (1 - 2 r^2) \left(\left(1 - r^2 \right) \alpha_{11} + r \sqrt{1 - r^2} \alpha_{31} - r \sqrt{1 - r^2} \alpha_{13} - r^2 \alpha_{33} \right) \\
    &= 2 r^2 \left(1 - r^2 \right) \alpha_{11} + 2 r^3 \sqrt{1 - r^2} \alpha_{31} + 2 r (1 - r^2) \sqrt{1 - r^2} \alpha_{13} + 2 r^2 (1 - r^2) \alpha_{33}. \\
    \implies& \left(r \sqrt{1 - r^2} c \phi_2 + r \sqrt{1 - r^2} \right) c \phi_3 - \sqrt{1 - r^2} s \phi_2 s \phi_3 \\
    &= r \sqrt{1 - r^2} \alpha_{11} + r^2 \alpha_{31} + \left(1 - r^2 \right) \alpha_{13} + r \sqrt{1 - r^2} \alpha_{33}.
\end{split}
\end{align}
Similar to the derivation of the closed-form expression for $\phi_1$ for the $LGR$ path, the above equation can be used to solve for $\phi_3$ for any case apart from $\phi_2 = \pi.$ In this case, the coefficient of both $c \phi_3$ and $s \phi_3$ is zero. This case is considered separately. For any other case, dividing both sides of Eq.~\eqref{eq: LGR_solving_phi_3} by $\sqrt{\left(r \sqrt{1 - r^2} c \phi_2 + r \sqrt{1 - r^2} \right)^2 + (1 - r^2) s^2 \phi_2},$
Eq.~\eqref{eq: LGR_solving_phi_3} can be rewritten as
\begin{align}
    c (\phi_3 + \beta) &= \frac{r \sqrt{1 - r^2} \alpha_{11} + r^2 \alpha_{31} + \left(1 - r^2 \right) \alpha_{13} + r \sqrt{1 - r^2} \alpha_{33}}{\sqrt{\left(r \sqrt{1 - r^2} c \phi_2 + r \sqrt{1 - r^2} \right)^2 + (1 - r^2) s^2 \phi_2}}. \\
    \therefore \phi_3 &= \left(\cos^{-1}{\left(\frac{r \sqrt{1 - r^2} \alpha_{11} + r^2 \alpha_{31} + \left(1 - r^2 \right) \alpha_{13} + r \sqrt{1 - r^2} \alpha_{33}}{\sqrt{\left(r \sqrt{1 - r^2} c \phi_2 + r \sqrt{1 - r^2} \right)^2 + (1 - r^2) s^2 \phi_2}} \right)} - \beta \right) \% 2 \pi,
\end{align}
where $\beta = \text{atan2} \left(s \phi_2, r \left(c \phi_2 + 1 \right) \right)$. From the above equation, for each solution of $\phi_2,$ at most two solutions can be obtained for $\phi_3$.

\subsubsection{Case with $\phi_2 = \pi$}

In this case, it was observed that $\phi_1$ and $\phi_3$ cannot be solved by pre- and post-multiplying by vectors. Hence, this case is addressed by solving the matrix equation given in Eq.~\eqref{eq: LGR_general_path_solve} directly, which is simplified as
\begin{align} \label{eq: LGR_rL_rR_equal_phi_2_pi}
    R_L (r, \phi_1) R_G (\pi) R_R (r, \phi_3) = \begin{pmatrix}
        \alpha_{11} & \alpha_{12} & \alpha_{13} \\
        \alpha_{21} & \alpha_{22} & \alpha_{23} \\
        \alpha_{31} & \alpha_{32} & \alpha_{33}
    \end{pmatrix}.
\end{align}
Expanding the LHS,
\begin{align*}
    &R_L (r, \phi_1) R_G (\pi) R_R (r, \phi_3) \\
    &= R_L (r, \phi_1) \begin{pmatrix}
        -1 & 0 & 0 \\
        0 & -1 & 0 \\
        0 & 0 & 1
    \end{pmatrix} \begin{pmatrix}
        1 - (1 - c \phi_3) r^2 & - r s \phi_3 & -(1 - c \phi_3) r \sqrt{1 - r^2} \\
        r s \phi_3 & c \phi_3 & s \phi_3 \sqrt{1 - r^2} \\
        -(1 - c \phi_3) r \sqrt{1 - r^2} & -s \phi_3 \sqrt{1 - r^2} & c \phi_3 + (1 - c \phi_3) r^2
    \end{pmatrix} \\
    &= \begin{pmatrix}
        1 - (1 - c \phi_1) r^2 & - r s \phi_1 & (1 - c \phi_1) r \sqrt{1 - r^2} \\
        r s \phi_1 & c \phi_1 & - s \phi_1 \sqrt{1 - r^2} \\
        (1 - c \phi_1) r \sqrt{1 - r^2} & s \phi_1 \sqrt{1 - r^2} & c \phi_1 + (1 - c \phi_1) r^2
    \end{pmatrix} \begin{pmatrix}
        -1 + (1 - c \phi_3) r^2 & r s \phi_3 & (1 - c \phi_3) r \sqrt{1 - r^2} \\
        -r s \phi_3 & -c \phi_3 & -s \phi_3 \sqrt{1 - r^2} \\
        -(1 - c \phi_3) r \sqrt{1 - r^2} & -s \phi_3 \sqrt{1 - r^2} & c \phi_3 + (1 - c \phi_3) r^2
    \end{pmatrix} \\
    &= \begin{pmatrix}
        -1 + r^2 - r^2 c (\phi_1 + \phi_3) & r s (\phi_1 + \phi_3) & r \sqrt{1 - r^2} (1 - c (\phi_1 + \phi_3)) \\
        -r s (\phi_1 + \phi_3) & -c (\phi_1 + \phi_3) & -\sqrt{1 - r^2} s (\phi_1 + \phi_3) \\
        -r \sqrt{1 - r^2} (1 - c (\phi_1 + \phi_3)) & -\sqrt{1 - r^2} s (\phi_1 + \phi_3) & r^2 + (1 - r^2) c (\phi_1 + \phi_3)
    \end{pmatrix}.
\end{align*}
Therefore,
\begin{align}
    \tan{\left(\phi_1 + \phi_3 \right)} = \frac{\alpha_{12}}{-r \alpha_{22}}.
\end{align}
Noting that infinitely many solutions exist, and the left and right turns have the same cost, $\phi_3$ can be set to zero without loss of generality, and $\phi_1$ can be solved from the above equation.

\newpage
\subsection{RGL path}

An $RGL$ path can be constructed using the derived $LGR$ path. To this end, assuming without loss of generality that the initial configuration corresponds to the identity matrix, we reflect the final configuration about the $XY$ plane, construct the $LGR$ path, and use the obtained parameters for the $RGL$ path by reflecting it about the $XY$ plane. To this end, the final configuration is modified to be
\begin{align} \label{eq: modified_final_configuration}
    R_{f, mod} = \begin{pmatrix}
        \alpha_{11} & \alpha_{12} & -\alpha_{13} \\
        \alpha_{21} & \alpha_{22} & -\alpha_{23} \\
        -\alpha_{31} & -\alpha_{32} & \alpha_{33}
    \end{pmatrix}.
\end{align}
An $LGR$ path is constructed to the modified final configuration. The obtained values for $\phi_1,$ $\phi_2,$ and $\phi_3$ are the parameters of the $RGL$ path to attain the initially provided final configuration (before modification).

\subsection{LRL path}

Consider an $LRL$ path, wherein the angle of the first $L$ segment is $\phi_1,$ the angle of the middle $R$ is equal to $\phi_2,$ and the angle of the final $L$ segment is $\phi_3$. The equation to be solved is given by
\begin{align} \label{eq: LRL_general_path_solve}
    R_L (r, \phi_1) R_R (r, \phi_2) R_L (r, \phi_3) = \begin{pmatrix}
        \alpha_{11} & \alpha_{12} & \alpha_{13} \\
        \alpha_{21} & \alpha_{22} & \alpha_{23} \\
        \alpha_{31} & \alpha_{32} & \alpha_{33}
    \end{pmatrix}.
\end{align}
Pre-multiplying Eq.~\eqref{eq: LRL_general_path_solve} by $\mathbf{u}_{L}^T$ and post-multiplying by $\mathbf{u}_{L},$
\begin{align} \label{eq: uLLTLRLuLL}
    \mathbf{u}_{L}^T R_L (r, \phi_1) R_R (r, \phi_2) R_L (r, \phi_3) \mathbf{u}_{L} = \mathbf{u}_{L}^T R_R (r, \phi_2) \mathbf{u}_{L} &= \mathbf{u}^T_{L} \begin{pmatrix}
        \alpha_{11} & \alpha_{12} & \alpha_{13} \\
        \alpha_{21} & \alpha_{22} & \alpha_{23} \\
        \alpha_{31} & \alpha_{32} & \alpha_{33}
    \end{pmatrix} \mathbf{u}_{L}.
\end{align}
Both sides of the above equation can be expanded to obtain
\begin{align*}
    \left(1 - 2 r^2 \right)^2 + 4 r^2 \left(1 - r^2 \right) \cos{\phi_2} = \left(1 - r^2 \right) \alpha_{11} + r \sqrt{1 - r^2} \left(\alpha_{13} + \alpha_{31} \right) + r^2 \alpha_{33}.
\end{align*}
The above equation can be rearranged to obtain
\begin{align*}
    \cos{\phi_2} = \frac{\left(1 - r^2 \right) \alpha_{11} + r \sqrt{1 - r^2} \left(\alpha_{13} + \alpha_{31} \right) + r^2 \alpha_{33} - \left(1 - 2 r^2 \right)^2}{4 r^2 \left(1 - r^2 \right)}.
\end{align*}
At most one solution for $\phi_2 \in (\pi, 2 \pi)$ can be obtained from the above equation.

Now, it is desired to obtain the solution for $\phi_1.$ Pre-multiplying Eq.~\eqref{eq: LRL_general_path_solve} by $\mathbf{u}_R^T$ and post-multiplying by $\mathbf{u}_L,$ the equation obtained is given by
\begin{align}
    \mathbf{u}_{R}^T R_L (r, \phi_1) R_R (r, \phi_2) R_L (r, \phi_3) \mathbf{u}_{L} = \mathbf{u}_{R}^T R_L (r, \phi_1) R_R (r, \phi_2) \mathbf{u}_{L} &= \mathbf{u}^T_{R} \begin{pmatrix}
        \alpha_{11} & \alpha_{12} & \alpha_{13} \\
        \alpha_{21} & \alpha_{22} & \alpha_{23} \\
        \alpha_{31} & \alpha_{32} & \alpha_{33}
    \end{pmatrix} \mathbf{u}_{L}.
\end{align}
Expanding both sides of the above equation, the equation obtained is given by
\begin{align*}
    &4 r^2 \left(1 - r^2 \right) \sin(\phi_1) \sin(\phi_2) + 4\left(2 r^2 - 1 \right) \left(1 - r^2 \right) r^2 \left(1 - \cos\left(\phi_2\right) \right) \cos(\phi_1) \\
    & + 8 r^6 - 12 r^4 + 6 r^2 - 1 - 4 \left(2 r^6 - 3 r^4 + r^2 \right) \cos(\phi_2) = \left(r^2 - 1 \right) \alpha_{11} + r \sqrt{1 - r^2} \left(\alpha_{31} - \alpha_{13} \right) + r^2 \alpha_{33}.
\end{align*}
The above equation can be rewritten as
\begin{align*}
    &\left(2 r^2 - 1 \right) \left(1 - \cos{\phi_2} \right) \cos{\phi_1} + \sin{\phi_2} \sin{\phi_1}\\
    &= \frac{\left(r^2 - 1 \right) \alpha_{11} + r \sqrt{1 - r^2} \left(\alpha_{31} - \alpha_{13} \right) + r^2 \alpha_{33} - \left(8 r^6 - 12 r^4 + 6 r^2 - 1 - 4 \left(2 r^6 - 3 r^4 + r^2 \right) \cos(\phi_2) \right)}{4 r^2 \left(1 - r^2 \right)}.
\end{align*}
Noting that $\sin{\phi_2} \neq 0,$ the coefficient of $\sin{\phi_1} \neq 0.$ Therefore, at most two solutions can be obtained for $\phi_1$ from the above equation for each $\phi_2$.

Now, it is desired to obtain the solutions for $\phi_3.$ To this end, pre-multiplying Eq.~\eqref{eq: LRL_general_path_solve} by $\mathbf{u}_{L}^T$ and post-multiplying by $\mathbf{u}_{R},$ the equation obtained is given by
\begin{align}
    \mathbf{u}_{L}^T R_L (r, \phi_1) R_R (r, \phi_2) R_L (r, \phi_3) \mathbf{u}_{R} = \mathbf{u}_{L}^T R_R (r, \phi_2) R_L (r, \phi_3) \mathbf{u}_{R} &= \mathbf{u}^T_{L} \begin{pmatrix}
        \alpha_{11} & \alpha_{12} & \alpha_{13} \\
        \alpha_{21} & \alpha_{22} & \alpha_{23} \\
        \alpha_{31} & \alpha_{32} & \alpha_{33}
    \end{pmatrix} \mathbf{u}_{R}.
\end{align}
Both sides of the above equation can be expanded to obtain
\begin{align*}
    &4 r^2 \left(1 - r^2 \right) \sin(\phi_3) \sin(\phi_2) + 4\left(2 r^2 - 1 \right) \left(1 - r^2 \right) r^2 \left(1 - \cos\left(\phi_2\right) \right) \cos(\phi_3) \\
    & + 8 r^6 - 12 r^4 + 6 r^2 - 1 - 4 \left(2 r^6 - 3 r^4 + r^2 \right) \cos(\phi_2) = \left(r^2 - 1 \right) \alpha_{11} + r \sqrt{1 - r^2} \left(\alpha_{13} - \alpha_{31} \right) + r^2 \alpha_{33}.
\end{align*}
The above equation can be rewritten as
\begin{align*}
    &\left(2 r^2 - 1 \right) \left(1 - \cos{\phi_2} \right) \cos{\phi_3} + \sin{\phi_2} \sin{\phi_3}\\
    &= \frac{\left(r^2 - 1 \right) \alpha_{11} + r \sqrt{1 - r^2} \left(\alpha_{13} - \alpha_{31} \right) + r^2 \alpha_{33} - \left(8 r^6 - 12 r^4 + 6 r^2 - 1 - 4 \left(2 r^6 - 3 r^4 + r^2 \right) \cos(\phi_2) \right)}{4 r^2 \left(1 - r^2 \right)}.
\end{align*}
Noting that $\sin{\phi_2} \neq 0,$ the coefficient of $\sin{\phi_3} \neq 0.$ Therefore, at most two solutions can be obtained for $\phi_1$ from the above equation for each $\phi_2$.

\subsection{RLR path}

An $RLR$ path is constructed by swapping the initial and final configurations, and swapping the tangent vector and tangent-normal directions. That is, the initial configuration and final configuration are modified to be
\begin{align*}
    R_i &= \begin{pmatrix}
        \alpha_{11} & -\alpha_{12} & -\alpha_{13} \\
        \alpha_{21} & -\alpha_{22} & -\alpha_{23} \\
        \alpha_{31} & -\alpha_{32} & -\alpha_{33}
    \end{pmatrix}, \\
    R_f &= \begin{pmatrix}
        1 & 0 & 0 \\
        0 & -1 & 0 \\
        0 & 0 & -1
    \end{pmatrix}.
\end{align*}
An $LRL$ path is then constructed. The obtained parameters $\phi_1, \phi_2,$ and $\phi_3$ are the angles of the final $R$ segment, the middle $L$ segment, and first $R$ segment of the $RLR$ path.

\subsection{$LR_\pi L$ path}

Consider an $LR_\pi L$ path, wherein the angle of the first $L$ segment is $\phi_1,$ the angle of the middle $R$ segment is $\pi,$ and the angle of the final $L$ segment is $\phi_3$. The equation to be solved is given by
\begin{align} \label{eq: LRpiL_general_path_solve}
    R_L (r, \phi_1) R_R (r, \pi) R_L (r, \phi_3) = \begin{pmatrix}
        \alpha_{11} & \alpha_{12} & \alpha_{13} \\
        \alpha_{21} & \alpha_{22} & \alpha_{23} \\
        \alpha_{31} & \alpha_{32} & \alpha_{33}
    \end{pmatrix}.
\end{align}
Pre-multiplying Eq.~\eqref{eq: LRpiL_general_path_solve} by $\mathbf{u}_{R}^T$ and post-multiplying by $\mathbf{u}_{L},$ the equation obtained is given by
\begin{align*}
    &\mathbf{u}_{R}^T R_L (r, \phi_1) R_R (r, \pi) R_L (r, \phi_3) \mathbf{u}_{L} \\
    &= \mathbf{u}_{R}^T R_L (r, \phi_1) R_R (r, \pi) \mathbf{u}_{L} = \mathbf{u}_{R}^T \begin{pmatrix}
        \alpha_{11} & \alpha_{12} & \alpha_{13} \\
        \alpha_{21} & \alpha_{22} & \alpha_{23} \\
        \alpha_{31} & \alpha_{32} & \alpha_{33}
    \end{pmatrix} \mathbf{u}_{L}.
\end{align*}
Expanding both sides of the above equation, the equation obtained in terms of $\phi_1$ is given by
\begin{align*}
    &- \left(2 r^2 - 1 \right) \left(-8 r^4 + 8 \left(r^2 - 1 \right) r^2 \cos{\phi_1} + 8 r^2 - 1 \right) = \alpha_{11} \left(r^2 - 1 \right) + r \left(\left(\alpha_{31} - \alpha_{13} \right) \sqrt{1 - r^2} + \alpha_{33} r \right).
\end{align*}
In the case of $r \neq \frac{1}{\sqrt{2}},$ which is a special case that will be handled separately, the above equation can be rewritten as
\begin{align}
    \cos{\phi_1} &= \frac{1}{8 \left(r^2 - 1 \right) r^2} \left(1 - 8 r^2 + 8 r^4 + \frac{\alpha_{11} \left(r^2 - 1 \right) + r \left(\left(\alpha_{31} - \alpha_{13} \right) \sqrt{1 - r^2} + \alpha_{33} r \right)}{1 - 2 r^2} \right).
\end{align}
The above equation yields at most two solutions for $\phi_1.$

The equation for $\phi_3$ can be obtained by pre-multiplying Eq.~\eqref{eq: LRpiL_general_path_solve} by $\mathbf{u}_{L}^T$ and post-multiplying by $\mathbf{u}_{R}$ as
\begin{align*}
    &\mathbf{u}_{L}^T R_L (r, \phi_1) R_R (r, \pi) R_L (r, \phi_3) \mathbf{u}_{R} \\
    &= \mathbf{u}_{L}^T R_R (r, \pi) R_L (r, \phi_3) \mathbf{u}_{L} = \mathbf{u}_{L}^T \begin{pmatrix}
        \alpha_{11} & \alpha_{12} & \alpha_{13} \\
        \alpha_{21} & \alpha_{22} & \alpha_{23} \\
        \alpha_{31} & \alpha_{32} & \alpha_{33}
    \end{pmatrix} \mathbf{u}_{R}.
\end{align*}
Expanding both sides of the above equation, the equation obtained in terms of $\phi_3$ is given by
\begin{align*}
    &- \left(2 r^2 - 1 \right) \left(-8 r^4 + 8 \left(r^2 - 1 \right) r^2 \cos{\phi_3} + 8 r^2 - 1 \right) = \alpha_{11} \left(r^2 - 1 \right) + r \left(\left(\alpha_{13} - \alpha_{31} \right) \sqrt{1 - r^2} + \alpha_{33} r \right).
\end{align*}
In the case of $r \neq \frac{1}{\sqrt{2}},$ which is a special case that will be handled separately, the above equation can be rewritten as
\begin{align}
    \cos{\phi_3} &= \frac{1}{8 \left(r^2 - 1 \right) r^2} \left(1 - 8 r^2 + 8 r^4 + \frac{\alpha_{11} \left(r^2 - 1 \right) + r \left(\left(\alpha_{13} - \alpha_{31} \right) \sqrt{1 - r^2} + \alpha_{33} r \right)}{1 - 2 r^2} \right).
\end{align}
The above equation yields at most two solutions for $\phi_3.$

\subsubsection{Special case of $r = \frac{1}{\sqrt{2}}$}

In this case, the net rotation matrix corresponding to the $LR_\pi L$ path reduces to
\begin{align*}
    R_L \left(\frac{1}{\sqrt{2}}, \phi_1 \right) R_R \left(\frac{1}{\sqrt{2}}, \pi \right) R_L \left(\frac{1}{\sqrt{2}}, \phi_3 \right) = \begin{pmatrix} 
        \frac{1}{2} (\cos{\left(\phi_1 - \phi_3 \right)} - 1) & \frac{\sin(\phi_1 - \phi_3)}{\sqrt{2}} & -\cos ^2\left(\frac{\phi_1-\phi_3}{2}\right) \\
        \frac{\sin (\phi_1-\phi_3)}{\sqrt{2}} & -\cos (\phi_1-\phi_3) & -\frac{\sin (\phi_1-\phi_3)}{\sqrt{2}} \\
        -\cos ^2\left(\frac{\phi_1-\phi_3}{2}\right) & -\frac{\sin (\phi_1-\phi_3)}{\sqrt{2}} & \frac{1}{2} (\cos{\left(\phi_1 - \phi_3 \right)} - 1).
    \end{pmatrix}
\end{align*}
Noting that this matrix must equal the RHS matrix, $\phi_1 - \phi_3$ can be obtained as
\begin{align*}
    \phi_1 - \phi_3 = \text{atan2} \left(\sqrt{2} \alpha_{21}, -\alpha_{22} \right),
\end{align*}
which yields an angle in $(-\pi, \pi].$ If the obtained angle is negative, then $\phi_1$ is set to zero and $\phi_3$ is computed. If the obtained angle is positive, $\phi_3$ is set to zero, and $\phi_1$ is computed.

\subsection{$RL_\pi R$ path}

The construction of an $RL_\pi R$ path using the construction of the $LR_\pi L$ follows similar to the construction of an $RLR$ path using the $LRL$ path construction. In particular, the final configuration is reflected about the $XY$ plane, the $LR_\pi L$ path is constructed, the parameters of which correspond to the $RL_\pi R$ path prior to reflection.

\subsection{LRLR path}

Consider an $LRLR$ path, wherein the angle of the first $L$ segment is $\phi_1,$ the angle of the middle $R$ and $L$ segments are equal to $\phi_2,$ and the angle of the final $R$ segment is $\phi_3$. The equation to be solved is given by
\begin{align} \label{eq: LRLR_general_path_solve}
    R_L (r, \phi_1) R_R (r, \phi_2) R_L (r, \phi_2) R_R (r, \phi_3) = \begin{pmatrix}
        \alpha_{11} & \alpha_{12} & \alpha_{13} \\
        \alpha_{21} & \alpha_{22} & \alpha_{23} \\
        \alpha_{31} & \alpha_{32} & \alpha_{33}
    \end{pmatrix}.
\end{align}
Pre-multiplying Eq.~\eqref{eq: LRLR_general_path_solve} by $\mathbf{u}_{L}^T$ and post-multiplying by $\mathbf{u}_{R},$
\begin{align*}
    &\mathbf{u}_{L}^T R_L (r, \phi_1) R_R (r, \phi_2) R_L (r, \phi_2) R_R (r, \phi_3) \mathbf{u}_{R} \\
    &= \mathbf{u}_{L}^T R_R (r, \phi_2) R_L (r, \phi_2) \mathbf{u}_{R} = \mathbf{u}_{L}^T \begin{pmatrix}
        \alpha_{11} & \alpha_{12} & \alpha_{13} \\
        \alpha_{21} & \alpha_{22} & \alpha_{23} \\
        \alpha_{31} & \alpha_{32} & \alpha_{33}
    \end{pmatrix} \mathbf{u}_{R}.
\end{align*}
Expanding both sides of the above equation, the equation obtained in terms of $\phi_2$ is given by
\begin{align*}
    &-1 + 10 r^2 - 16 r^4 + 8 r^6 -8 \left(r^2-3 r^4+2 r^6 \right) \cos{\phi_2} + 8 r^4 \left(r^2 - 1 \right) \cos^2{\phi_2} \\
    &= \alpha_{11} \left(r^2 - 1\right) + r \left(\alpha_{13} \sqrt{1 - r^2} - \alpha_{31} \sqrt{1 - r^2} + \alpha_{33} r \right).
\end{align*}
Noting that the above equation is a quadratic equation in terms of $\cos{\phi_2},$ at most two real solutions can be obtained from $\cos{\phi_2}.$ Furthermore, noting that for every value of $\cos{\phi_2},$ two solutions exist, wherein one solution lies in $[0, \pi],$ and another solution lies in $[\pi, 2 \pi),$ at most one solution for $\phi_2$ is selected since $\phi_2 \in (\pi, 2 \pi)$ for optimality. Therefore, from the obtained equation, at most two solutions for $\phi_2 \in (\pi, 2 \pi)$ can be obtained.

Consider pre-multiplying Eq.~\eqref{eq: LRLR_general_path_solve} by $\mathbf{u}_R^T$ and post-multiplying by $\mathbf{u}_R.$ The equation obtained is given by
\begin{align*}
    &\mathbf{u}_{R}^T R_L (r, \phi_1) R_R (r, \phi_2) R_L (r, \phi_2) R_R (r, \phi_3) \mathbf{u}_{R} \\
    &= \mathbf{u}_{R}^T R_L (r, \phi_1) R_R (r, \phi_2) R_L (r, \phi_2) \mathbf{u}_{R} = \mathbf{u}_{R}^T \begin{pmatrix}
        \alpha_{11} & \alpha_{12} & \alpha_{13} \\
        \alpha_{21} & \alpha_{22} & \alpha_{23} \\
        \alpha_{31} & \alpha_{32} & \alpha_{33}
    \end{pmatrix} \mathbf{u}_{R}.
\end{align*}
Expanding both sides of the above equation, the equation obtained in terms of $\phi_1$ is given by
\begin{align*}
    &\left(2 r^2 - 1\right) \left(12 r^6 - 20 r^4 + 10 r^2 + 4 \left(r^2 - 1\right) r^4 \cos{(2 \phi_2)} - 8 \left(2 r^6 - 3 r^4 + r^2 \right) \cos{(\phi_2)} - 1 \right) \\
    & -4 r^2 \left(r^2 - 1\right) \left(2 r^2 \cos(\phi_2) - 2 r^2 + 1 \right) \left(\cos(\phi_1) \left(\left(2 r^2 - 1 \right) \cos(\phi_2) - 2 r^2 + 2 \right) - \sin(\phi_1) \sin(\phi_2) \right) \\
    &= \alpha_{11} \left(1 - r^2 \right) + r \left(\alpha_{13} \left(-\sqrt{1 - r^2} \right) - \alpha_{31} \sqrt{1 - r^2} + \alpha_{33} r \right).
\end{align*}
The above equation can be rewritten as
\begin{align} \label{eq: solving_phi1_LRLR_path}
    A \cos{\phi_1} + B \sin{\phi_1} + C = \alpha_{11} \left(1 - r^2 \right) + r \left(\alpha_{13} \left(-\sqrt{1 - r^2} \right) - \alpha_{31} \sqrt{1 - r^2} + \alpha_{33} r \right),
\end{align}
where
\begin{align*}
    A &= 4 r^2 \left(1 - r^2 \right) \left(2 r^2 \cos(\phi_2) - 2 r^2 + 1 \right) \left(\left(2 r^2 - 1 \right) \cos(\phi_2) - 2 r^2 + 2 \right), \\
    B &= 4 r^2 \left(1 - r^2 \right) \left(2 r^2 \cos(\phi_2) - 2 r^2 + 1 \right) \left(-\sin{\phi_2} \right),\\
    C &= \left(2 r^2 - 1\right) \left(12 r^6 - 20 r^4 + 10 r^2 + 4 \left(r^2 - 1\right) r^4 \cos{(2 \phi_2)} - 8 \left(2 r^6 - 3 r^4 + r^2 \right) \cos{(\phi_2)} - 1 \right).
\end{align*}
\textbf{In the case that $2 r^2 \cos{\phi_2} - 2 r^2 + 1 = 0,$ $\phi_1$ cannot be uniquely solved from the above equation. This special case will be considered separately.}

Suppose $2 r^2 \cos{\phi_2} - 2 r^2 + 1 \neq 0.$ It is desired to determine if $A^2 + B^2 \neq 0$ to obtain a finite number of solutions from the above equation. 
It should be noted that if $A^2 + B^2 = 0,$ then $A = 0, B = 0.$ However, for $B = 0,$ it is necessary that $\sin{\phi_2} = 0.$ However, for an $LRLR$ path to be optimal, $\phi_2 \in (\pi, 2\pi).$ Hence, since $B \neq 0,$ $A^2 + B^2 \neq 0.$ Therefore, at most two solutions can be obtained for $\phi_1$ from Eq.~\eqref{eq: solving_phi1_LRLR_path}.

Now, it is desired to obtain the expression for $\phi_3.$ Pre-multiplying Eq.~\eqref{eq: LRLR_general_path_solve} by $\mathbf{u}_L^T$ and post-multiplying by $\mathbf{u}_L,$ the equation obtained is given by
\begin{align*}
    &\mathbf{u}_{L}^T R_L (r, \phi_1) R_R (r, \phi_2) R_L (r, \phi_2) R_R (r, \phi_3) \mathbf{u}_{L} \\
    &= \mathbf{u}_{L}^T R_R (r, \phi_2) R_L (r, \phi_2) R_R (r, \phi_3) \mathbf{u}_{L} = \mathbf{u}_{L}^T \begin{pmatrix}
        \alpha_{11} & \alpha_{12} & \alpha_{13} \\
        \alpha_{21} & \alpha_{22} & \alpha_{23} \\
        \alpha_{31} & \alpha_{32} & \alpha_{33}
    \end{pmatrix} \mathbf{u}_{L}.
\end{align*}
Expanding both sides of the above equation, the equation obtained in terms of $\phi_3$ is given by
\begin{align*}
    A \cos{\phi_3} + B \sin{\phi_3} + C = \left(1 - r^2 \right) \alpha_{11} + r \sqrt{1 - r^2} \left(\alpha_{13} + \alpha_{31} \right) + r^2 \alpha_{33},
\end{align*}
where the expressions for $A, B,$ and $C$ are the same as obtained for the equation for $\phi_1.$ At most two solutions can be obtained when $2 r^2 \cos{\phi_2} - 2 r^2 + 1 \neq 0$ for each $\phi_2.$

\subsubsection{Special case when $\cos{\phi_2} = 1 - \frac{1}{2 r^2}$}

In this case, $\phi_1$ and $\phi_3$ cannot be uniquely solved. Noting that using $\cos{\phi_2} = 1 - \frac{1}{2 r^2}$, a solution for $\phi_2 \in (\pi, 2 \pi)$ exists for $r \in \left(\frac{1}{2}, 1 \right)$ it is desired to determine whether the path exists or not. Since $\phi_2 \in (\pi, 2 \pi),$ the corresponding value of $\sin{\phi_2}$ can be obtained as $\sin{\phi_2} = -\sqrt{1 - \cos^2{\phi_2}} = \frac{-\sqrt{4 r^2 - 1}}{2 r^2}.$ Using these values of $\sin{\phi_2}$ and $\cos{\phi_2},$ the net rotation matrix in Eq.~\eqref{eq: LRLR_general_path_solve} reduces to
\begin{align*}
    R_L (r, \phi_1) R_R (r, \phi_2) R_L (r, \phi_2) R_R (r, \phi_3) = \begin{pmatrix}
        \gamma_{11} & \gamma_{12} & \gamma_{13} \\
        -\gamma_{12} & \gamma_{22} & \gamma_{23} \\
        -\gamma_{13} & \gamma_{23} & \gamma_{33}
    \end{pmatrix},
\end{align*}
where
\begin{align*}
    \gamma_{11} &= \frac{1}{2} \left(\sqrt{4 r^2-1} \sin (\phi_1+\phi_3)+\left(1-2 r^2\right) \cos (\phi_1+\phi_3)+2 r^2-2\right), \\
    \gamma_{12} &= \frac{\left(2 r^2-1\right) \sin (\phi_1+\phi_3)+\sqrt{4 r^2-1} \cos (\phi_1+\phi_3)}{2 r}, \\
    \gamma_{13} &= \frac{\sqrt{1-r^2} \left(\sqrt{4 r^2-1} \sin (\phi_1+\phi_3)+\left(1-2 r^2\right) \cos (\phi_1+\phi_3)+2 r^2\right)}{2 r}, \\
    \gamma_{22} &= \frac{\sqrt{4 r^2-1} \sin (\phi_1+\phi_3)+\left(1-2 r^2\right) \cos (\phi_1+\phi_3)}{2 r^2}, \\
    \gamma_{23} &= -\frac{\sqrt{1-r^2} \left(\left(2 r^2-1\right) \sin (\phi_1+\phi_3)+\sqrt{4 r^2-1} \cos (\phi_1+\phi_3)\right)}{2 r^2}, \\
    \gamma_{33} &= \frac{2 r^4+\left(r^2-1\right) \sqrt{4 r^2-1} \sin (\phi_1+\phi_3)+\left(-2 r^4+3 r^2-1\right) \cos (\phi_1+\phi_3)}{2 r^2}.
\end{align*}
% \begin{align*}
%     \begin{pmatrix}
%         \frac{1}{2} \left(\sqrt{4 r^2-1} \sin (\phi_1+\phi_3)+\left(1-2 r^2\right) \cos (\phi_1+\phi_3)+2 r^2-2\right) & \frac{\left(2 r^2-1\right) \sin (\phi_1+\phi_3)+\sqrt{4 r^2-1} \cos (\phi_1+\phi_3)}{2 r} & \frac{\sqrt{1-r^2} \left(\sqrt{4 r^2-1} \sin (\phi_1+\phi_3)+\left(1-2 r^2\right) \cos (\phi_1+\phi_3)+2 r^2\right)}{2 r} \\ \frac{\left(1-2 r^2\right) \sin (\phi_1+\phi_3)-\sqrt{4 r^2-1} \cos (\phi_1+\phi_3)}{2 r} & \frac{\sqrt{4 r^2-1} \sin (\phi_1+\phi_3)+\left(1-2 r^2\right) \cos (\phi_1+\phi_3)}{2 r^2} & -\frac{\sqrt{1-r^2} \left(\left(2 r^2-1\right) \sin (\phi_1+\phi_3)+\sqrt{4 r^2-1} \cos (\phi_1+\phi_3)\right)}{2 r^2} \\ \frac{\sqrt{1-r^2} \left(-\sqrt{4 r^2-1} \sin (\phi_1+\phi_3)+\left(2 r^2-1\right) \cos (\phi_1+\phi_3)-2 r^2\right)}{2 r} & -\frac{\sqrt{1-r^2} \left(\left(2 r^2-1\right) \sin (\phi_1+\phi_3)+\sqrt{4 r^2-1} \cos (\phi_1+\phi_3)\right)}{2 r^2} & \frac{2 r^4+\left(r^2-1\right) \sqrt{4 r^2-1} \sin (\phi_1+\phi_3)+\left(-2 r^4+3 r^2-1\right) \cos (\phi_1+\phi_3)}{2 r^2}
%     \end{pmatrix}
% \end{align*}
Clearly, either no solution or infinitely many solutions can be obtained for $\phi_1$ and $\phi_3,$ since unique solutions cannot be obtained for $\phi_1$ and $\phi_3$. 

If infinitely many solutions are obtained, then, $\phi_3$ can be set to zero since the left and right turns have the same radius and cost. Then, using $\gamma_{12}$ and $\gamma_{22},$
\begin{align*}
    \begin{pmatrix}
        \frac{\sqrt{4 r^2 - 1}}{2 r} & \frac{2 r^2 - 1}{2 r} \\
        -\frac{2 r^2 - 1}{2 r^2} & \frac{\sqrt{4 r^2 - 1}}{2 r^2} 
    \end{pmatrix} \begin{pmatrix}
        \cos{\left(\phi_1 \right)} \\
        \sin{\left(\phi_1 \right)}
    \end{pmatrix} = \begin{pmatrix}
        \alpha_{12} \\
        \alpha_{22}
    \end{pmatrix}.
\end{align*}
The determinant of the matrix in the left-hand side is $r^2,$ which is non-zero. Hence, the expression for $\cos{\phi_1}$ and $\sin{\phi_1}$ can be obtained as
\begin{align*}
    \begin{pmatrix}
        \cos{\left(\phi_1 \right)} \\
        \sin{\left(\phi_1 \right)}
    \end{pmatrix} = \frac{1}{r^2} \begin{pmatrix}
        \frac{\sqrt{4 r^2 - 1}}{2 r^2} & -\frac{2 r^2 - 1}{2 r} \\
        \frac{2 r^2 - 1}{2 r^2} & \frac{\sqrt{4 r^2 - 1}}{2 r} 
    \end{pmatrix} \begin{pmatrix}
        \alpha_{12} \\
        \alpha_{22}
    \end{pmatrix}.
\end{align*}
It should be noted that the obtained expressions for $\cos{\phi_1}$ and $\sin{\phi_1}$ must be verified to satisfy $\cos^2{\phi_1} + \sin^2{\phi_1} = 1.$

\subsection{RLRL path}

The construction of an $RLRL$ path can be performed using the construction of the $LRLR$ path by first reflecting the final configuration about the $XY$ plane (assuming the initial configuration is the identity matrix without loss of generality). The $LRLR$ path can then be constructed to attain the modified final configuration. The obtained solutions for $\phi_1,$ $\phi_2,$ and $\phi_3$ will correspond to the arc angles of the first $R$ segment, intermediary $L$ and $R$ segments, and final $L$ segment in the $RLRL$ path connecting to the initially provided final configuration. To this end, the final configuration is modified as given in Eq.~\eqref{eq: modified_final_configuration}, using which the $LRLR$ path is constructed.

\subsection{LRLRL path}

Consider an $LRLRL$ path, wherein the angle of the first $L$ segment is $\phi_1,$ the angle of the middle $R, L,$ and $R$ segments are equal to $\phi_2,$ and the angle of the final $L$ segment is $\phi_3$. The equation to be solved is given by
\begin{align} \label{eq: LRLRL_general_path_solve}
    R_L (r, \phi_1) R_R (r, \phi_2) R_L (r, \phi_2) R_R (r, \phi_2) R_L (r, \phi_3) = \begin{pmatrix}
        \alpha_{11} & \alpha_{12} & \alpha_{13} \\
        \alpha_{21} & \alpha_{22} & \alpha_{23} \\
        \alpha_{31} & \alpha_{32} & \alpha_{33}
    \end{pmatrix}.
\end{align}
Pre-multiplying Eq.~\eqref{eq: LRLRL_general_path_solve} by $\mathbf{u}_{L}^T$ and post-multiplying by $\mathbf{u}_{L},$
\begin{align*}
    &\mathbf{u}_{L}^T R_L (r, \phi_1) R_R (r, \phi_2) R_L (r, \phi_2) R_R (r, \phi_2) R_L (r, \phi_3) \mathbf{u}_{L} \\
    &= \mathbf{u}_{L}^T R_R (r, \phi_2) R_L (r, \phi_2) R_R (r, \phi_2) \mathbf{u}_{L} = \mathbf{u}_{L}^T \begin{pmatrix}
        \alpha_{11} & \alpha_{12} & \alpha_{13} \\
        \alpha_{21} & \alpha_{22} & \alpha_{23} \\
        \alpha_{31} & \alpha_{32} & \alpha_{33}
    \end{pmatrix} \mathbf{u}_{L}.
\end{align*}
Expanding both sides of the above equation, the equation obtained in terms of $\phi_2$ is given by
\begin{align*}
    &16 r^8-48 r^6+48 r^4-16 r^2+1 - 16 r^2 \left(1 - r^2 \right)^2 \left(3 r^2 - 1 \right) \cos{\phi_2} + 16 r^4 \left(2 - 5 r^2 + 3 r^4 \right) \cos^2{\phi_2} \\
    &+ 16 r^6 \left(1 - r^2 \right) \cos^3{\phi_2} = \alpha_{11} \left(1 - r^2 \right) + r \left(\alpha_{13} \sqrt{1 - r^2} + \alpha_{31} \sqrt{1 - r^2} + \alpha_{33} r\right).
\end{align*}

At most three real solutions can be obtained for $\cos{\phi_2}.$ The solutions that lie in $[-1, 1]$ are selected; for each solution of $\cos{\phi_2},$ at most one solution lies in $(\pi, 2 \pi).$ It is now desired to obtain the corresponding solutions for $\phi_1$ and $\phi_3.$

Pre-multiplying Eq.~\eqref{eq: LRLRL_general_path_solve} by $\mathbf{u}_R^T$ and post-multiplying by $\mathbf{u}_L,$
\begin{align*}
    &\mathbf{u}_{R}^T R_L (r, \phi_1) R_R (r, \phi_2) R_L (r, \phi_2) R_R (r, \phi_2) R_L (r, \phi_3) \mathbf{u}_{L} \\
    &= \mathbf{u}_{R}^T R_L (r, \phi_1) R_R (r, \phi_2) R_L (r, \phi_2) R_R (r, \phi_2) \mathbf{u}_{L} = \mathbf{u}_{R}^T \begin{pmatrix}
        \alpha_{11} & \alpha_{12} & \alpha_{13} \\
        \alpha_{21} & \alpha_{22} & \alpha_{23} \\
        \alpha_{31} & \alpha_{32} & \alpha_{33}
    \end{pmatrix} \mathbf{u}_{L}.
\end{align*}
Expanding both sides of the above equation, the equation obtained in terms of $\phi_1$ is given by
\begin{align*}
    A \cos{\phi_1} + B \sin{\phi_1} + C = \left(r^2 - 1 \right) \alpha_{11} + r \sqrt{1 - r^2} \left(\alpha_{31} - \alpha_{13} \right) + r^2 \alpha_{33},
\end{align*}
where
\begin{align*}
    A &= \gamma \tan\left(\frac{\phi_2}{2}\right)\left(-4r^4+2\left(2r^2-1\right)r^2\cos(\phi_2)+6r^2-1\right), \\
    B &= \gamma \left(1 - 2 r^2 + 2 r^2 \cos{\phi_2} \right), \\
    C &= \left(1 - 2 r^2 \right) \bigg[4 r^8 \cos(3\phi_2) - 40 r^8 -4 r^6 \cos(3\phi_2) + 88 r^6 - 64 r^4 + 16 r^2 - 8 \left(3 r^4 - 5 r^2 + 2 \right) r^4 \cos(2 \phi_2) \\
    & \hspace{2.1cm} + 4 \left(15 r^6 - 31 r^4 + 20 r^2 - 4 \right) r^2 \cos(\phi_2) - 1 \bigg], \\
    \gamma &= 8 r^2 \left(1 - r^2 \right) \left(1 - r^2 + r^2 \cos{\phi_2} \right) \sin{\phi_2}.
\end{align*}

It is claimed that $A^2 + B^2$ can be zero only when $\gamma = 0.$ When $\gamma = 0,$ the corresponding value of $\phi_2$ must correspond to $\cos{\phi_2} = 1 - \frac{1}{r^2}.$ This is a case that will be considered separately.

Consider $\gamma \neq 0.$ Then, $B = 0$ only when $\cos{\phi_2} = 1 - \frac{1}{2 r^2},$ which has a solution in $(\pi, 2 \pi)$ when $r \in \left(\frac{1}{2}, 1 \right].$ Since $\phi_2 \in (\pi, 2 \pi),$ the corresponding expression for $\sin{\phi_2} = -\sqrt{1 - \cos^2{\phi_2}} = \frac{-\sqrt{4 r^2 - 1}}{2 r^2}.$ Since $\tan{\left(\frac{\phi_2}{2} \right)} \neq 0,$ $A$ can be zero only if $-4 r^4 + 2 \left(2 r^2 - 1 \right) r^2 \cos{\phi_2} + 6 r^2 - 1$ reduces to zero.
Substituting the expression for $\cos{\phi_2}$ in the considered expression, we get
\begin{align*}
    -4 r^4 + 2 \left(2 r^2 - 1 \right) r^2 \cos{\phi_2} + 6 r^2 - 1 = 2 r^2,
\end{align*}
which is non-zero. Hence, $A^2 + B^2 \neq 0$ if $\gamma \neq 0.$ Therefore, at most two solutions can be obtained for $\phi_1$ for each solution of $\phi_2.$

\textbf{Remark:} For the implementation, the expressions for $A$ and $B$ can be alternately written as
\begin{align*}
    A &= 16 r^2 \left(r^2 - 1 \right) \sin^2\left(\frac{\phi_2}{2} \right) \left(-6 r^6 + 11 r^4 - 7 r^2 + \left(r^4 - 2 r^6 \right) \cos(2\phi_2) + \left(8 r^4 - 12 r^2 + 3 \right) r^2 \cos(\phi_2) + 1 \right), \\
    B &= 8r^2 \left(1 - r^2 \right) \sin(\phi_2) \left(r^4 \cos(2\phi_2) + 3 r^4 - 3 r^2 + \left(3 r^2 - 4 r^4 \right) \cos(\phi_2) + 1 \right).
\end{align*}

It is now desired to obtain the expression for $\phi_3.$ Pre-multiplying Eq.~\eqref{eq: LRLRL_general_path_solve} by $\mathbf{u}_L^T$ and post-multiplying by $\mathbf{u}_R,$
\begin{align*}
    &\mathbf{u}_{L}^T R_L (r, \phi_1) R_R (r, \phi_2) R_L (r, \phi_2) R_R (r, \phi_2) R_L (r, \phi_3) \mathbf{u}_{R} \\
    &= \mathbf{u}_{L}^T R_R (r, \phi_2) R_L (r, \phi_2) R_R (r, \phi_2) R_L (r, \phi_3) \mathbf{u}_{R} = \mathbf{u}_{L}^T \begin{pmatrix}
        \alpha_{11} & \alpha_{12} & \alpha_{13} \\
        \alpha_{21} & \alpha_{22} & \alpha_{23} \\
        \alpha_{31} & \alpha_{32} & \alpha_{33}
    \end{pmatrix} \mathbf{u}_{R}.
\end{align*}
Expanding both sides of the above equation, the equation obtained in terms of $\phi_3$ is given by
\begin{align*}
    A \cos{\phi_1} + B \sin{\phi_1} + C = \left(r^2 - 1 \right) \alpha_{11} + r \sqrt{1 - r^2} \left(\alpha_{13} - \alpha_{31} \right) + r^2 \alpha_{33},
\end{align*}
where the expressions for $A, B,$ and $C$ are the same as that obtained for the equation for $\phi_1.$ Hence, at most two solutions can be obtained for $\phi_3$ for each $\phi_2$ if $\gamma \neq 0.$

\subsubsection{Case with $\cos{\phi_2} = 1 - \frac{1}{r^2}$}

In this case, the net rotation matrix corresponding to the $LRLRL$ path, given by the LHS of Eq.~\eqref{eq: LRLRL_general_path_solve}, reduces to
\begin{align*}
    R_L (r, \phi_1) R_R (r, \phi_2) R_L (r, \phi_2) R_R (r, \phi_2) R_L (r, \phi_3) = \begin{pmatrix}
        \delta_{11} & \delta_{12} & \delta_{13} \\
        -\delta_{12} & \delta_{22} & \delta_{23} \\
        \delta_{13} & -\delta_{23} & \delta_{33}
    \end{pmatrix},
\end{align*}
where $\sin{\phi_2} = -\sqrt{1 - \cos^2{\phi_2}} = -\frac{\sqrt{2r^2 - 1}}{r^2}$ was used. Here,
\begin{align*}
    \delta_{11} &= -2\sin\left(\frac{\phi_1+\phi_3}{2}\right)\left(\left(r^2-1\right)\sin\left(\frac{\phi_1+\phi_3}{2}\right)+\sqrt{2r^2-1}\cos\left(\frac{\phi_1+\phi_3}{2}\right)\right), \\
    \delta_{12} &= -\frac{\left(r^2-1\right)\sin(\phi_1+\phi_3)+\sqrt{2r^2-1}\cos(\phi_1+\phi_3)}{r}, \\
    \delta_{13} &= \frac{\sqrt{1-r^2}\left(\sqrt{2r^2-1}\sin(\phi_1+\phi_3)-\left(r^2-1\right)\cos(\phi_1+\phi_3)+r^2\right)}{r}, \\
    \delta_{22} &= \frac{\left(r^2-1\right)\cos(\phi_1+\phi_3)-\sqrt{2r^2-1}\sin(\phi_1+\phi_3)}{r^2}, \\
    \delta_{23} &= -\frac{\sqrt{1-r^2}\left(\left(r^2-1\right)\sin(\phi_1+\phi_3)+\sqrt{2r^2-1}\cos(\phi_1+\phi_3)\right)}{r^2}, \\
    \delta_{33} &= \frac{r^4+\left(r^2-1\right)\sqrt{2r^2-1}\sin(\phi_1+\phi_3)-\left(r^2-1\right)^2\cos(\phi_1+\phi_3)}{r^2}.
\end{align*}

Noting that the net rotation matrix must equal the desired final configuration, given by
\begin{align*}
    R_f = \begin{pmatrix}
        \alpha_{11} & \alpha_{12} & \alpha_{13} \\
        \alpha_{21} & \alpha_{22} & \alpha_{23} \\
        \alpha_{31} & \alpha_{32} & \alpha_{33}
    \end{pmatrix},
\end{align*}
the expression for $\phi_1 + \phi_3$ can be obtained by equation $\delta_{12}$ to $\delta_{22}$ to $\alpha_{12}$ and $\alpha_{22},$ respectively, to obtain
\begin{align*}
    \begin{pmatrix}
        \sqrt{2 r^2 - 1} & \left(r^2 - 1 \right) \\
        \left(r^2 - 1 \right) & -\sqrt{2 r^2 - 1}
    \end{pmatrix} \begin{pmatrix}
        \cos{\left(\phi_1 + \phi_3 \right)} \\
        \sin{\left(\phi_1 + \phi_3 \right)}
    \end{pmatrix} = \begin{pmatrix}
        -r \alpha_{12} \\
        r^2 \alpha_{22}
    \end{pmatrix}.
\end{align*}
The determinant of the matrix on the left hand side is $-r^4,$ which is non-zero. Hence, the matrix on the left hand side can be inverted to obtain
\begin{align*}
    \begin{pmatrix}
        \cos{\left(\phi_1 + \phi_3 \right)} \\
        \sin{\left(\phi_1 + \phi_3 \right)}
    \end{pmatrix} &= \frac{-1}{r^4} \begin{pmatrix}
        -\sqrt{2 r^2 - 1} & -\left(r^2 - 1 \right) \\
        -\left(r^2 - 1 \right) & \sqrt{2 r^2 - 1}
    \end{pmatrix} \begin{pmatrix}
        -r \alpha_{12} \\
        r^2 \alpha_{22}
    \end{pmatrix} \\
    &= \frac{-1}{r^4} \begin{pmatrix}
        \sqrt{2 r^2 - 1} r \alpha_{12} - r^2 \left(r^2 - 1 \right) \alpha_{22} \\
        r \left(r^2 - 1 \right) \alpha_{12} + r^2 \sqrt{2 r^2 - 1} \alpha_{22}
    \end{pmatrix}.
\end{align*}
Setting $\phi_3 = 0$ without loss of generality, a solution for $\phi_1$ can be obtained from the above equation. It should be noted that the obtained expressions for $\cos{\phi_1}$ and $\sin{\phi_1}$ must be verified to satisfy $\cos^2{\phi_1} + \sin^2{\phi_1} = 1.$

\subsection{RLRLR path}

The $RLRLR$ path can be constructed by utilizing the construction of an $LRLRL$ path. To this end, the final configuration is reflected about the $XY$ plane, and the $LRLRL$ path is constructed to the reflected final configuration. The modified final configuration is given in Eq.~\eqref{eq: modified_final_configuration}. The obtained solutions for $\phi_1,$ $\phi_2,$ and $\phi_3$ correspond to the parameters for the $RLRLR$ path.

\bibliographystyle{IEEEtran}
\bibliography{IEEEabrv, references}

\begin{thebibliography}{1}
\providecommand{\url}[1]{#1}
\csname url@samestyle\endcsname
\providecommand{\newblock}{\relax}
\providecommand{\bibinfo}[2]{#2}
\providecommand{\BIBentrySTDinterwordspacing}{\spaceskip=0pt\relax}
\providecommand{\BIBentryALTinterwordstretchfactor}{4}
\providecommand{\BIBentryALTinterwordspacing}{\spaceskip=\fontdimen2\font plus
\BIBentryALTinterwordstretchfactor\fontdimen3\font minus
  \fontdimen4\font\relax}
\providecommand{\BIBforeignlanguage}[2]{{%
\expandafter\ifx\csname l@#1\endcsname\relax
\typeout{** WARNING: IEEEtran.bst: No hyphenation pattern has been}%
\typeout{** loaded for the language `#1'. Using the pattern for}%
\typeout{** the default language instead.}%
\else
\language=\csname l@#1\endcsname
\fi
#2}}
\providecommand{\BIBdecl}{\relax}
\BIBdecl

\bibitem{kumar2025newapproachmotionplanning}
\BIBentryALTinterwordspacing
D.~P. Kumar, S.~Darbha, S.~G. Manyam, and D.~Casbeer, ``A new approach to
  motion planning in 3d for a dubins vehicle: Special case on a sphere,'' 2025.
  [Online]. Available: \url{https://arxiv.org/abs/2504.01215}
\BIBentrySTDinterwordspacing

\end{thebibliography}

\end{document}